\documentclass[12pt]{article}
\usepackage{amsmath,amssymb,amstext,dsfont,fancyvrb,float,fontenc,graphicx,caption,subcaption,theorem,hyperref}
\usepackage[utf8]{inputenc}

\usepackage{booktabs, algorithm, algorithmic, url, hyperref}

\usepackage[letterpaper]{geometry}
\ifx\volno\undefined\def\volno{9}\fi
\ifx\volyear\undefined\def\volyear{2026}\fi
\ifx\pagno\undefined\def\pagno{242--274}\fi

\newfont{\footsc}{cmcsc10 at 8truept}
\newfont{\footbf}{cmbx10 at 8truept}
\newfont{\footrm}{cmr10 at 10truept}

\usepackage{fancyhdr}
\usepackage{relsize}
\usepackage{sectsty}
\allsectionsfont{\larger[-1]} 

\renewcommand\paragraph{\@startsection{paragraph}{4}{\z@}
                                    {2ex \@plus.5ex \@minus.2ex}
                                    {-1em}
                                    {\normalfont\normalsize\bfseries}}

\renewcommand\subparagraph{\@startsection{subparagraph}{5}{\parindent}
                                       {2ex \@plus.5ex \@minus .2ex}
                                       {-1em}
                                      {\normalfont\normalsize\bfseries}}

\newlength{\BiblioSpacing}
\renewenvironment{thebibliography}[1]{
\begin{oldthebibliography}{#1}
\setlength{\parskip}{\BiblioSpacing}
\setlength{\itemsep}{\BiblioSpacing}
}
{
\end{oldthebibliography}
}

\usepackage[strict]{changepage}
\def\abstractname{Abstract -}   
\def\abstract{\begin{adjustwidth}{1cm}{1cm} \par    \footnotesize \noindent {\bf \abstractname} 
\def\endabstract{ \end{adjustwidth} \smallskip }}

{\theorembodyfont{\itshape}\newtheorem{theorem}{Theorem}[section]}
{\theorembodyfont{\itshape}}
{\theorembodyfont{\itshape}}
{\theorembodyfont{\itshape}\newtheorem{lemma}[theorem]{Lemma}}
{\theorembodyfont{\itshape}}
{\theorembodyfont{\rm}}
{\theorembodyfont{\rm}}
{\theorembodyfont{\rm}}
{\theorembodyfont{\rm }}

\def\dedicatory{\date}

\newcommand{\qed}{\hfill$\square$}

\title{\Large\bf Generalizations of the M\textup{\&}M Game}
\author{\sc S. Das, 
S.J. Miller, 
G. Polanco, \\  
\sc Y. Wu, 
X. Wang, 
A. Yang, 
C. Yao\thanks{This work began during the 2024 Polymath Junior REU program and was supported in part by NSF Grant DMS2341670. We thank the participants at Polymath Junior, our colleagues, and the anonymous referee for helpful comments on this paper.}
}
\dedicatory{}

\begin{document}
\setcounter{page}{1}
\maketitle
\thispagestyle{fancy}

\vskip 1.5em

\begin{abstract}
The M\&M Game involves two players, each starting with a pile of M\&M's. During each round, each player tosses a fair coin: if the coin lands heads, that player eats one M\&M, and if it lands tails, the player does not eat. If, at the end of a round, one player still has M\&M's while the other has none, then the player with M\&M's remaining is declared the winner. If both players eat their last M\&M in the same round, the game is said to end in a tie. Previous work on the M\&M Game studied the probability of a tie and derived a simple closed-form expression in the special case where both players start with the same number of M\&M's. We generalize the M\&M Game in several directions, including allowing players to toss two coins per round and modifying the probability distributions of the coin flips. We use generating functions, Monte Carlo methods, and nonlinear curve fitting to study the generalized M\&M Game. 
\end{abstract}
 

\begin{MSC}
60C05
\end{MSC}

\section{Introduction} 
\label{sec:introduction} 
\subsection{Background of the M\&M Game} 
\label{subsection: background of the MM Game}
The M\&M Game arose from attempts to model life expectancy, inspired by a question asked by Cameron and Kayla Miller (aged 4 and 2): If two people are born on the same day, do they die on the same day? To explore this question, their father, Steven Miller, who happened to have bags of M\&M's with him, reformulated it as a two-player game. \emph{Suppose two players start with the same number of M\&M's. On each turn, each player tosses a coin and eats one M\&M if the coin lands heads and eats none if it lands tails. What is the probability that they eat their last M\&M's at the same time}? Later, they and Badinski, Huffaker, McCue, and Stone~\cite{badinski2017m} explored several different approaches to analyze the game, which turned out to be a great springboard for exploring many powerful mathematical techniques: hypergeometric functions, memoryless processes, and walks on graphs, to name a few. These connections inspired the first paper and this sequel, which introduces generating functions to analyze generalizations. 

From \cite{badinski2017m}, the rules of the M\&M Game are as follows. Player A and Player B begin with $I_1$ and $I_2$ M\&M's, respectively, both positive integers. During each round, each player flips a fair coin independently and simultaneously.\footnote{We retain the assumption that both players flip fair coins simultaneously only to remain faithful to the original description of the M\&M Game in \cite{badinski2017m}. In fact, the order of the flips is irrelevant; all that matters is that the players toss their coins independently.} If a player's coin lands heads, that player eats one M\&M; if it lands tails, that player does not eat an M\&M. The game continues until at least one player has no M\&M's remaining. If, at the end of a round, one player still has M\&M's while the other has none, then the player with M\&M's remaining is declared the winner. If both players run out of M\&M's in the same round, the game ends in a tie. In \cite{badinski2017m}, the authors study the probability of a tie when both players start with $I$ M\&M's, and we briefly summarize their work below.

Let $\mathbb{P}(I_1, I_2)$ denote the probability of a tie when Player A starts with $I_1$ M\&M's and Player B starts with $I_2$ M\&M's. In \cite{badinski2017m}, the authors first studied the probability of a tie when both players start with the same number of M\&M's and obtained
\begin{equation}
    \label{eqs: definition of P(k, k)}
    \mathbb{P}(I, I) \ =\ \sum_{k=I}^{\infty} \binom{k - 1}{I - 1}^2 \left(\frac{1}{2}\right)^{2k}.
\end{equation}
Their argument is as follows.  First, they decompose the tie probability according to the round on which the game ends: 
\begin{equation}
\label{eqs: decompose P(k, k) into P_n(k, k)}
\mathbb{P}(I, I) \ = \ \sum_{k= I}^{\infty} \mathbb{P}_k(I, I), 
\end{equation}
where $\mathbb{P}_k(I, I)$ denotes the probability that the game ends in a tie after exactly $k$ rounds, given that both players start with $I$ M\&M's. The sum starts at $k = I$ because \emph{at least} $I$ rounds are required for the game to end in a tie. It remains to compute $\mathbb{P}_k(I,I)$. For the game to end in a tie on the $k^{\text{th}}$ round, each player must eat exactly $I-1$ M\&M's during the first $k-1$ rounds and then eat her last M\&M on the $k^{\text{th}}$ round. As each player tosses only one coin in each round, this means that each player must obtain $I - 1$ heads in the first $k - 1$ tosses, followed by a head on the $k^{\text{th}}$ toss. Thus, there are $\binom{k - 1}{I - 1}\binom{1}{1}$ ways to accomplish this. Since the coin is fair and the two players flip independently, $\mathbb{P}_{k}(I, I)$ is the product of the probabilities that each player obtains the $I^{\text{th}}$ head at the $k^{\text{th}}$ toss, and this is given by 
\begin{equation} 
\label{eqs: definition of P_n(k, k)}
\mathbb{P}_{k}(I, I) \ = \ \binom{k - 1}{I - 1}\left( \frac{1}{2} \right)^{k} \binom{k - 1}{I - 1}\left( \frac{1}{2} \right)^{k} .
\end{equation}
Substituting Equation \eqref{eqs: definition of P_n(k, k)} into Equation \eqref{eqs: decompose P(k, k) into P_n(k, k)} gives Equation \eqref{eqs: definition of P(k, k)}, as desired.

Although Equation \eqref{eqs: definition of P(k, k)} gives the probability of a tie for the game, the formula is somewhat unsatisfactory: it is expressed as an infinite sum, which is difficult to evaluate, and it is hard to see how the probability of a tie depends on the starting number of M\&M's. By adopting a different perspective, they derived a sum of finitely many terms to compute the probability of a tie. The key insight in \cite{badinski2017m} is that the M\&M Game is a \textit{memoryless process}: the probability of a tie starting at a given state only depends on how many M\&M's each player currently possesses, not on how many rounds it takes for them to get there. By exploiting this memoryless property, they obtained the following result.  
\begin{theorem}[Page 201 of \cite{badinski2017m}] \label{thm: probability tie finite sum}  
The probability that the M\&M Game ends in a tie when two players each start with $I$ M\&M's is
$$
\mathbb{P}(I, I) \ = \ \sum_{k = 0}^{I - 1} \binom{2I - k - 2}{k} \binom{2I-2k-2}{I- k -1} \left(\frac{1}{3}\right)^{2I - k - 1}.
$$
\end{theorem}

To gain further insight into the game, it is helpful to view the M\&M Game as a random walk on the two-dimensional integer lattice. Formally, we represent the state of the game at each round by the lattice point $(m,n)$, where $m$ is the number of M\&M's currently held by Player A and $n$ is the number held by Player B. We use the notation $m$ and $n$ to denote the numbers of M\&M's currently possessed by the two players in a given round. This notation allows us to distinguish these quantities from $I_1$ and $I_2$, which denote the numbers of M\&M's with which the two players begin. As we shall see later, this distinction proves useful in our study of the generalized M\&M Game. From the perspective of random walks, each possible outcome of a round in the M\&M Game corresponds to one of the following moves on the lattice. 
\begin{itemize}
    \item \emph{Both Players Eat}: send $(m,n)$ to $(m-1,n-1)$.
    \item \emph{Only Player A Eats}: send $(m, n)$ to $(m-1,n)$.
    \item \emph{Only Player B Eats}: send $(m,n)$ to $(m,n-1)$.
    \item \emph{No Player Eats}: leave $(m,n)$ unchanged.
\end{itemize}
Since each player independently eats one M\&M with probability $1/2$ and does not eat with probability $1/2$, each of these four lattice moves occurs with probability $1/4$.

Thus, the tie probability in the M\&M Game can be reformulated as follows: starting from $(I_1,I_2)$, what is the probability of reaching the origin $(0,0)$ using the four moves described above? By the memoryless property of the M\&M Game, the sequence of moves that led to the current state $(m,n)$ is irrelevant when calculating the tie probability. Moreover, since the tie probability depends only on the number of M\&M's currently possessed by the two players, the number of null moves, corresponding to the outcome \emph{No Player Eats}, does not affect the tie probability. Thus, we may remove this move and rescale the probabilities of the remaining moves, which simplifies the problem. After removing the null move, the remaining three moves are equally likely with probability $1/3$. If we let $F(m,n)$ denote the probability that the walk eventually reaches $(0,0)$ when the current position is $(m,n)$, then $F(m,n)$ satisfies the recurrence relation\footnote{Even without removing the move \emph{No Player Eats} via the memoryless argument, we can still arrive at the recurrence relation \eqref{eqs: recurrence definition for the original MM Game}. Given that we are currently at $(m, n)$ and have four equally likely moves, we have the recurrence  
\[
F(m, n) \ =  \frac{1}{4} F(m - 1, n) + \frac{1}{4}F(m, n - 1) + \frac{1}{4}F(m -1, n -1) + \frac{1}{4}F(m, n), 
\]
and performing simple algebra gives the recurrence relation \eqref{eqs: recurrence definition for the original MM Game}. Thus, we can arrive at the desired recurrence relation even without using the memoryless nature of the M\&M Game.}
\begin{equation}
\label{eqs: recurrence definition for the original MM Game}
F(m,n) \ = \ \frac{1}{3}F(m-1,n) + \frac{1}{3}F(m,n-1) + \frac{1}{3}F(m-1,n-1).
\end{equation}
In this recurrence, we observe $F(0,0)=1$, since the probability of reaching $(0,0)$ when we are already at $(0,0)$ is 1. We also observe $F(m,0)=F(0,n)=0$ for all positive integers $m$ and $n$, since once one player has no M\&M's remaining while the other still has some left, the game can no longer end in a tie.

\subsection{Generalizing the M\&M Game} 
\label{subsection: generalizing the MM Game}
The goal of this paper is to generalize the M\&M Game. In the original paper  \cite{badinski2017m}, each player flips a single fair coin at each turn. To generalize the game, we consider the following two methods for each player.  
\begin{itemize}

    \item \textbf{Method 1.} The player tosses 1 fair coin in each round. If it lands heads, then the player eats $d$ M\&M's; if it lands tails, then the player eats none.

    \item \textbf{Method 2.} The player tosses 2 fair coins in each round. If the first coin lands heads, then the player eats $d_1$ M\&M's; if the first coin lands tails, then the player eats none. If the second coin lands heads, then the player eats $d_2$ M\&M's; if it lands tails, then the player eats none. 
\end{itemize}

In both methods, if the player is ever supposed to eat more than she has left, then the player eats all the remaining M\&M's. We focus on the M\&M Game for two players with arbitrary starting values $I_1$ and $I_2$.  Using these two methods, we create three possible games for two players.
\begin{itemize}
\item \textbf{Game 1.} Both players play Method 1; that is, a player eats $d$ M\&M's if her coin lands heads and eats none if her coin lands tails. 

\item \textbf{Game 2.} Player A plays Method 1 and Player B plays Method 2. Thus, Player A eats $d$ M\&M's if her coin lands heads, while Player B eats $d_1$ M\&M's if her first coin lands heads and eats $d_2$ M\&M's if her second coin lands heads.

\item \textbf{Game 3.} Both players play Method 2; that is, a player eats $d_1$ M\&M's if her first coin lands heads and eats $d_2$ M\&M's if her second coin lands heads.
\end{itemize}
Following the rules of the M\&M Game in \cite{badinski2017m}, we say that the game ends in a tie if both players exhaust all of their M\&M's in the same round. We formalize this using the lattice representation. Let $(m,n)$ denote the numbers of M\&M's currently possessed by Player A and Player B, respectively. Suppose that in a given round Player A is supposed to eat $a$ M\&M's and Player B is supposed to eat $b$ M\&M's. If $m \leq a$ and $n \leq b$, then we say the game ends in a tie. If, on the other hand, one player exhausts all of her M\&M's while the other still has some M\&M's remaining, then we declare the player with M\&M's remaining to be the winner.

We can now ask several natural questions. What is the tie probability of each game? How does the tie probability depend on $d$, $d_1$, and $d_2$? The game becomes far more interesting and complicated now. To illustrate, we consider Game 2, in which Player A eats $d = 5$ M\&M's whenever a head is flipped, while Player B eats $d_1 = 2$ M\&M's if the first coin lands heads and eats $d_2 = 3$ M\&M's if the second coin lands heads. Suppose both players start with $4$ M\&M's. Then even though the average number of M\&M's removed on each turn is the same (namely,  2.5) for both methods, we see that Player A goes out in one turn with probability $1/2$, while Player B does so only with probability $1/4$. However, if both players start with 2 M\&M's, then while Player A still goes out in one turn with probability $1/2$, now Player B does so with probability $3/4$. We see that the two methods for playing the game influence the outcome of the game in a nontrivial way, a complication that did not appear in \cite{badinski2017m}. 

The rest of this paper is organized as follows. In Section~\ref{sec: Comparison of Three Games}, we analyze the three proposed games and address the questions raised above. In Section \ref{subsec: Problem Approach}, we introduce the technique of  generating functions and provide an overview of our method for analyzing the tie probabilities in the three proposed games. In Section \ref{subsection: results}, we present explicit formulas for the tie probability of each game. We then give a derivation of the explicit formula for the tie probability in Game 1 in Section \ref{subsection: derivation of game 1}; since the arguments for Games 2 and 3 are similar, we refer the reader to Appendix \ref{appendix: Proof of Results for Game 2 and Game 3} for the detailed proofs. Finally, in Section \ref{subsec:simulations for three games}, we present numerical results illustrating how the number of coins a player may toss affects the tie probabilities.

In Section \ref{sec: games under different probabilities}, we study how the tie probability changes when the coin-flip distribution is altered. We begin by giving a straightforward generalization of Theorem \ref{thm: probability tie finite sum} to the case of a biased coin that lands heads with probability $p$, and we examine how the tie probability depends on $p$. This extension assumes that coin flips are independent across rounds, with a fixed probability of landing heads. It is also natural to ask how the game changes when this probability depends on the history of previous rounds. For this reason, we introduce a state-dependent version of the M\&M Game in which the probability that a coin lands heads depends on the number of M\&M's remaining. We then study this game through extensive Monte Carlo simulations and find that the  tie probabilities are well approximated by a Gompertz curve.

Finally, in Section \ref{sec: conclusion}, we propose several directions for future work, including extensions to games with more complicated rules, the study of the expected number of turns, and more complex interactions between coin outcomes and the starting number of M\&M's.


\section{Comparison of Three Games}\label{sec: Comparison of Three Games}

\subsection{Problem Approach}
\label{subsec: Problem Approach}
\subsubsection{Recurrence Formulation for the M\&M Game}
We study the tie probabilities for the three games introduced in Section \ref{subsection: generalizing the MM Game}. Since both players begin with finitely many M\&M's, and since each player can consume only a nonnegative number of M\&M's on each turn, the game terminates after finitely many turns with probability one.\footnote{This relies on the assumption that each player has a positive probability of consuming at least one M\&M in each round; if this probability were $0$ for both players, the game could remain in the same state forever. Under the game rules considered in this section, both players have a positive probability of consuming some M\&M's. } This observation allows us to avoid the infinite-sum approach used in Equation \eqref{eqs: definition of P(k, k)} and instead formulate the problem using recurrence relations. As in \cite{badinski2017m}, we formulate the tie probability in terms of a recurrence relation, though the analysis here is technically more involved because the number
of M\&M’s removed in a turn may vary.

We first introduce some notation. Let $F(m,n)$ denote the probability of a tie when Player A and Player B currently have $m$ and $n$ M\&M's, respectively. In addition, let $d_i$ denote the number of M\&M's a player eats when her $i^{\text{th}}$ coin lands heads. When a player flips only one coin, we use $d$ to denote the number of M\&M's eaten when the coin lands heads. We represent the possible moves available to a player in a single round using multisets, which are like sets except that repeated elements are allowed. We denote multisets using square brackets $[\ ]$.\footnote{One might wonder why we use multisets rather than sets to describe the moves available to a player. The reason is that sets may fail to distinguish between different ways of obtaining the same move. For example, suppose a player uses Method 2 with $d_1=d_2=1$. Then the possible amounts she can eat are $0$, $d_1=1$, $d_2=1$, and $d_1+d_2=2$. If we used a set, these moves would be recorded as ${0,1,2}$, losing one occurrence of the move $1$. In contrast, using a multiset records them as $[0,1,1,2]$, preserving the fact that the move $1$ can occur in two distinct ways.
} As described earlier, each player follows either Method 1 or Method 2. Let $\mathcal{A}$ denote the multiset of all possible moves available to Player A under her chosen method. If Player A uses Method 1, then on a given round Player A either eats $d$ M\&M's or none, depending on the outcome of Player A's coin flip; thus $\mathcal{A} = [0,d]$. If Player A plays Method 2, then Player A could eat $d_1$ M\&M's, $d_2$ M\&M's, $d_1 + d_2$ M\&M's, or no M\&M's, depending on the outcomes of Player A's two coin flips; thus, $\mathcal{A} = [0,d_1,d_2,d_1+d_2]$. We define the multiset $\mathcal{B}$ analogously for Player B.

Let $\mathcal{S}$ denote the multiset of all possible moves in a round after removing the null move $(0, 0)$. Formally,
$\mathcal{S} = (\mathcal{A} \times \mathcal{B}) \setminus [(0,0)]$,
where the null move $(0,0)$ is omitted because, as discussed in Section \ref{subsection: background of the MM Game}, it leaves the state unchanged and therefore does not affect the tie probability. From the perspective of random walks on the integer lattice, $\mathcal{S}$ may be interpreted as the multiset of admissible steps from one lattice point to another.

Finally, let $P(a,b)$ denote the probability that Player A eats $a$ M\&M's and Player B eats $b$ M\&M's in a single round. In the three games considered in this section, we may take $P(a,b) = 1/\left|\mathcal{S}\right|$, though this notation is also useful for possible future extensions in which different moves are assigned different probabilities.

By conditioning on the next move from the current state $(m, n)$, we can express the tie probability of each game as a recurrence relation of the form
\begin{equation}
F(m, n)  \ = \
\begin{cases}
\sum_{(a,b)\in\mathcal S} P(a,b)\,F(m-a,n-b)
& \text{if } m\ge 1 \text{ and } n\ge 1  \\[6pt]

0
& \text{if exactly one of } m,n \text{ is } \le 0  \\[6pt]

1
& \text{if } m\le 0 \text{ and } n\le 0,
\end{cases}
\label{eqs: general recurrence structure}
\end{equation}
where the latter two cases are the initial conditions. We now explain each initial condition. If exactly one of $m$ and $n$ is less than or equal to $0$, then one player has already exhausted her M\&M's while the other has not, so a tie is impossible. Thus, $F(m,n)=0$ if exactly one of $m, n \leq 0$. If $m=n=0$, then both players have exhausted their M\&M's in the same round, and the game ends in a tie. In our version of the game, it is possible for both players to consume more M\&M's than they have remaining, and we still treat such a case as a tie. Thus, $F(m,n)=1$ whenever $m \leq 0$ and $n \leq 0$.

\subsubsection{Sketch of the Argument}
We now sketch our argument and present the detailed derivation of the tie probabilities for Game 1 in Section \ref{subsection: derivation of game 1}. The arguments for Games 2 and 3 follow a similar approach and are given in Appendix \ref{appendix: Proof of Results for Game 2 and Game 3}.

For each game, we first define a multiset of moves $\mathcal{S}$ based on the methods used by the two players. Substituting $\mathcal{S}$ into Equation \eqref{eqs: general recurrence structure} yields the corresponding recurrence for the game under consideration. We then associate this recurrence with the generating function  $A(x, y) = \sum_{m, n \geq 0} F(m, n)x^my^n$. 

Our goal is to use the technique of generating functions to derive an explicit sum of finitely many terms for $F(m,n)$; see \cite{wilf_generatingfunctionology} for background on generating functions. Accordingly, we write $[x^m y^n]A(x,y)$ for the coefficient of $x^m y^n$ in $A(x,y)$. We then apply the following results to simplify $A(x,y)$.

\begin{lemma}
Under the recurrence relation \eqref{eqs: general recurrence structure}, we have
\begin{equation*}
A(x, y) \ = \ 1 + \sum_{m, n \geq 1} F(m, n)x^m y^n.
\end{equation*}
\label{lemma: A(x,x) = 1 + sum F(m, n)x^m y^n}
\end{lemma}

\begin{lemma}
Let $a, b \geq 1$. Under the recurrence relation \eqref{eqs: general recurrence structure}, we have
\begin{equation*}
\sum_{m, n \geq 1} F(m - a, n - b)x^m y^n \ = \  \sum^{a}_{M = 1}\sum^{b}_{N = 1} x^M y^N + x^{a}y^{b} \left( A(x, y) - 1\right).
\end{equation*}
Also, we have
\begin{equation*}
\sum_{m \geq 1} \sum_{n \geq 1} F(m - a, n)x^m y^n  \ = \ x^a \left(A(x, y) - 1\right)
\end{equation*}
and
\begin{equation*}
\sum_{m \geq 1} \sum_{n \geq 1} F(m, n - b)x^m y^n  \ = \ y^b \left(A(x, y) - 1\right).
\end{equation*}
\label{lemma: rewrite F(m -a, n - b)}
\end{lemma}
The proofs for Lemmas \ref{lemma: A(x,x) = 1 + sum F(m, n)x^m y^n} and \ref{lemma: rewrite F(m -a, n - b)} are algebraic manipulations using the recurrence relation and initial conditions. The detailed proofs can be found in Appendix \ref{appendix: Proof of Lemmas}.

The technique of generating functions and the two lemmas introduced above allow us to obtain another quick proof of Theorem \ref{thm: probability tie finite sum}. As a warm-up to the complicated derivations for Games 1, 2, and 3, we shall apply our argument to the original M\&M Game introduced in Section \ref{sec:introduction}. We observe that the step multiset is $\mathcal{S} = [(1, 0), (0, 1), (1, 1)]$, and inserting $\mathcal{S}$ into Equation \eqref{eqs: general recurrence structure} yields 
\begin{equation*}
F(m,n) \  =  \
\begin{cases}
\begin{aligned}[t]
\frac{1}{3}\Bigl[F(m-1,n)&+F(m,n-1)\\
& +F(m-1,n-1)\Bigr] 
\end{aligned}
& \text{if } m\ge 1 \text{ and } n\ge 1  \\[6pt]
0
&   
 \text{if exactly one of } m,n 
 \text{ is }\leq 0, \\[6pt]
1 & \text{if } m\le 0 \text{ and } n\le 0. 
\end{cases}
\end{equation*}
Associating this recurrence relation with $A(x, y) = \sum_{m, n \geq 0} F(m, n)x^my^n$, we obtain
\begin{align*}
A(x, y) &\ = \  1 + \sum_{m, n \geq 1} F(m, n)x^m y^n  \\  
&\ = \ 1 + \frac{1}{3}\biggl[ \sum_{m \geq 1}\sum_{n \geq 1}F(m - 1, n)x^m y^n   \nonumber \\
&\qquad \qquad + \sum_{m \geq 1} \sum_{n \geq 1}F(m, n - 1)x^m y^n + \sum_{m \geq 1} \sum_{n \geq 1}F(m - 1, n- 1)x^m y^n \biggr] \\ 
&\ = \ 1 + \frac{1}{3}\left[x\left(A(x, y) - 1\right) + y \left(A(x, y) - 1 \right) + xy + xy  \left(A(x, y) - 1\right) \right],
\end{align*}
where we applied Lemma \ref{lemma: A(x,x) = 1 + sum F(m, n)x^m y^n} in the first step and Lemma \ref{lemma: rewrite F(m -a, n - b)} in the last step. Solving for $A(x,y)$, we find  
\begin{equation}
\label{eqs: Delannoy Numbers} 
A(x,y) \ = \ 1 + \frac{\frac{1}{3}xy}{1 - \frac{1}{3}\left( x + y + xy\right)}. 
\end{equation}

We now extract the coefficients of $A(x, y)$. After a simple application of the geometric series formula and multinomial expansion,  we find 
\begin{align*}
A(x, y) &\ = \ 1+ xy \sum^{\infty}_{k = 0}\left(\frac{1}{3}\right)^{k + 1} \left(x + y + xy \right)^k \\
&\ = \ 1 +  \sum^{\infty}_{k = 0} \left(\frac{1}{3}\right)^{k + 1} \sum_{\substack{i, j, \ell \geq 0 \\ i + j + \ell = k}} \frac{k!}{i!j!\ell!} x^{i + \ell + 1} y^{j + \ell + 1}. 
\end{align*} 
We let $I_1 = i + \ell + 1$ and $I_2 = j + \ell + 1$, so we have
$$
A(x, y) \ = \  1 + \sum_{I_1, I_2 \geq 1} x^{I_1} y^{I_2}  \sum_{\ell = 0}^{\min(I_1 - 1, I_2 - 1)} \binom{I_1 + I_2 - \ell - 2}{\ell, I_1 - \ell - 1, I_2 - \ell - 1}\left(\frac{1}{3}\right)^{I_1 + I_2 -  \ell - 1}, 
$$
where we note that 
$$
\binom{I_1 + I_2 - \ell - 2}{\ell, I_1 - \ell - 1, I_2 - \ell - 1} \ = \ \binom{I_1 + I_2 - \ell - 2}{\ell} \binom{I_1 + I_2 - 2\ell -2}{I_1 - \ell - 1}.
$$
Thus, we see that the coefficient of $x^{I_1}y^{I_2}$ is given by 
\begin{equation}
[x^{I_1} y^{I_2}]A(x, y) \ = \  \sum^{\min(I_1 -1, I_2 - 1)}_{\ell = 0}\binom{I_1 + I_2 - \ell - 2}{\ell} \binom{I_1 + I_2 - 2 \ell - 2}{I_1 - \ell - 1}\left(\frac{1}{3}\right)^{I_1 + I_2  - \ell - 1}, 
\label{eqs: coefficients of the original MM Game}
\end{equation}
and specializing to $I_1 = I_2 = I$ gives the result stated in Theorem \ref{thm: probability tie finite sum}. 

The overall procedure for deriving the tie probabilities of Games 1, 2, and 3 is the same as the procedure illustrated above, but the computation is more involved. We describe the procedure here. We first insert the recurrence relation prescribed by the step multiset $\mathcal{S}$ into  the generating function $A(x, y) = \sum_{m, n \geq 0} F(m, n)x^m y^n$. Using Lemmas \ref{lemma: A(x,x) = 1 + sum F(m, n)x^m y^n} and \ref{lemma: rewrite F(m -a, n - b)}, we can rewrite the generating function as 
\begin{equation}
\label{eqs: generating function complicated expression}
    A(x, y) \ = \  1 + t B(\alpha, \beta, t) \sum_{(D_1, D_2) \in \mathcal{D}} \sum_{M = 1}^{D_1}\sum_{N = 1}^{D_2}x^M y^N.
\end{equation}
We now explain Equation \eqref{eqs: generating function complicated expression}. The variable $t$ is a constant and denotes the probability of each move in a single round; in our game setup, $t$ is $1/|\mathcal{S}|$. The terms $\alpha$ and $\beta$ are polynomials in $x$ and $y$, respectively, that encode the positive numbers of M\&M's a player can eat in a single round. For example, if Player A uses Method 2, then the positive numbers of M\&M's she can eat are $d_1$, $d_2$, and $d_1+d_2$, so $\alpha = x^{d_1} + x^{d_2} + x^{d_1 + d_2}$. 

The multiset $\mathcal{D}$ is a sub-multiset of $\mathcal{S}$, and it consists of moves $(a,b)$ in $\mathcal{S}$ with $a>0$ and $b>0$.
For example, if we consider Game 1, then $\mathcal{D} = [(d, d)]$. If we consider Game 3, then $\mathcal{D} = [d_1, d_2, d_1 + d_2] \times [d_1, d_2,  d_1 + d_2]$.  In addition, the term $B(\alpha, \beta, t)$ in Equation \eqref{eqs: generating function complicated expression} is defined by
\begin{equation}
B(\alpha, \beta, t) \ = \  \frac{1}{1- t\left(\alpha + \beta + \alpha \beta\right)}.
\end{equation} 
This is precisely the generating function for the Delannoy numbers (for further discussion of the Delannoy numbers, see \cite{Banderier_2005} and \url{https://oeis.org/A008288}). In our previous derivation, we see that the generating function for the Delannoy numbers is already present as a factor in the second term on the right hand side of Equation \eqref{eqs: Delannoy Numbers}. More generally, for constant weight $t$, the Delannoy numbers are defined by  
\begin{equation} 
D\left(m, n, t\right) \ =  \ \sum^{\min(m, n)}_{\ell = 0} \binom{m + n - \ell}{m - \ell, n - \ell, \ell} t^{m + n - \ell}. 
\end{equation}
Thus, the right hand side of Equation \eqref{eqs: coefficients of the original MM Game} can be written as $\frac{1}{3}  D(I_1-1,I_2-1,\frac{1}{3})$. As we shall see later, Delannoy numbers play an important role in our study of the generalized M\&M Game.

The next step is to make a change of variables so that every term in \eqref{eqs: generating function complicated expression} is expressed in terms of $x$ and $y$. As $B(\alpha, \beta, t)$ contains $\alpha$ and $\beta$, we have to make a change of variables to rewrite $B(\alpha, \beta, t)$ in terms of $x$ and $y$. The last step is simply to extract 
the coefficients corresponding to $x^{I_1}y^{I_2}$, which give us the desired result.

\subsection{Results}
\label{subsection: results} 
We first introduce some notation. Let $\mathbb{P}^{(i)}$ denote the tie probability for Game $i$. Let $\mathbb{P}_{(d)}^{(1)}(I_1,I_2)$ denote the tie probability in Game 1, where the players start with $I_1$ and $I_2$ M\&M's, respectively, and a player eats $d$ M\&M's upon flipping a head. Let $\mathbb{P}_{(d,d_1,d_2)}^{(2)}(I_1,I_2)$ denote the tie probability in Game 2, where Player A starts with $I_1$ M\&M's and eats $d$ M\&M's upon flipping a head, while Player B starts with $I_2$ M\&M's and eats $d_1$ M\&M's if her first coin lands heads and $d_2$ M\&M's if her second coin lands heads. Finally, let $\mathbb{P}_{(d_1,d_2)}^{(3)}(I_1,I_2)$ denote the tie probability in Game 3, where both players start with $I_1$ and $I_2$ M\&M's, respectively, and a player eats $d_1$ M\&M's if her first coin lands heads and $d_2$ M\&M's if her second coin lands heads. We have the following results. 

\begin{theorem} 
\label{thm: finite sum for Game 1}
\textbf{\textup{(Tie Probability for Game 1)}}
In Game 1, Player A and Player B start with $I_1$ and $I_2$ M\&M's, respectively. In each round, each player tosses one fair coin, and a player eats $d$ M\&M's if that player's coin lands heads and eats none if it lands tails. The probability of a tie in Game 1 is
\begin{equation*}
\mathbb{P}_{(d)}^{(1)}(I_1, I_2) \ = \  \frac{1}{3}D\left(\left\lfloor \frac{I_1 - 1}{d} \right\rfloor, \left\lfloor \frac{I_2 - 1}{d} \right\rfloor, \frac{1}{3} \right).  
\end{equation*} 
 \end{theorem}

Intuitively, we expect the tie probability in Game 1 to agree with that in the original M\&M Game after scaling the initial number of M\&M's by a factor of $d$. For example, consider Game 1 with $d = 3$. If both players start with $6$ M\&M's in Game 1, then this game has the same tie probability as the original M\&M Game in which both players start with $2$ M\&M's. In Theorem \ref{thm: connection between Game 1 and original MM Game}, we make this connection precise.
\begin{theorem}
In the original M\&M Game, Player A and Player B both start with $I$ M\&M's. In each round, each player tosses one fair coin, and a player eats one M\&M  if that player's coin lands heads  and eats none if it lands tails. Let $\mathbb{P}^{\mathrm{original}}\left(I, I\right)$ denote the tie probability in the original M\&M Game when both players start with $I$ M\&M's. We then have 
\begin{equation*} 
\mathbb{P}^{\mathrm{original}}\left(I, I\right) \ = \ \mathbb{P}_{(d)}^{(1)}(dI, dI). 
\end{equation*} 
\label{thm: connection between Game 1 and original MM Game}
\end{theorem}

\begin{theorem} 
\textbf{\textup{(Tie Probability for Game 2)}}
In Game 2, Player A and Player B start with $I_1$ and $I_2$ M\&M's, respectively. In each round, Player A tosses one fair coin and eats $d$ M\&M's if her coin lands heads. Player B tosses two fair coins, eating $d_1$ M\&M's if her first coin lands heads and $d_2$ M\&M's if her second coin lands heads. The probability of a tie in Game 2 is
\begin{align*}
\mathbb{P}_{(d, d_1, d_2)}^{(2)}(I_1, I_2) &\ = \ \frac{1}{7} \Biggl( \sum_{(D_1, D_2) \in \mathcal{D}} \sum_{\substack{i, j, k \geq 0 \\ 1 \leq I_2 - d_1(i + k) - d_2(j + k) \leq D_2}}\binom{i+ j + k}{i, j, k} 
\\ 
&\qquad \qquad \qquad \cdot D\left(\left\lfloor \frac{I_1 - 1}{d} \right\rfloor, i + j + k, \frac{1}{7} \right) \Biggr),   
\end{align*}
where $\mathcal{D} = \left[d \right] \times \left[ d_1, d_2, d_1 + d_2\right]$.
\label{thm: finite sum for Game 2}
\end{theorem}

\begin{theorem} 
\textbf{\textup{(Tie Probability for Game 3)}}
In Game 3, Player A and Player B start with $I_1$ and $I_2$ M\&M's, respectively. In each round, each player tosses two fair coins. If the player's first coin lands heads, that player eats $d_1$ M\&M's, and if the player's second coin lands heads, that player eats $d_2$ M\&M's. The probability of a tie in Game 3 is
\begin{align*}
\mathbb{P}_{(d_1, d_2)}^{(3)}(I_1, I_2)&\ = \ \frac{1}{15}\Biggl(\sum_{(D_1, D_2) \in \mathcal{D}}  \sum_{\substack{r, s, t \geq 0 \\ 1 \leq I_1 - d_1(r + t) - d_2(s + t) \leq D_1}}   \sum_{\substack{i, j, k \geq 0 \\ 1 \leq I_2 - d_1(i + k) - d_2(j + k) \leq D_2}}  \binom{r + s + t}{r,s,t} \nonumber \\
& \qquad \qquad \cdot  \binom{i + j + k}{i, j, k}   D\left(r + s + t, i + j + k, \frac{1}{15}\right)\Biggl),
\end{align*}
where $\mathcal{D} = \left[d_1, d_2, d_{1} + d_{2}\right] \times \left[d_1, d_2, d_1 + d_2\right]$. 
\label{thm: finite sum for Game 3}
\end{theorem}

\begin{figure}[!htbp]
    \centering   
    \includegraphics[width=0.7\textwidth ]{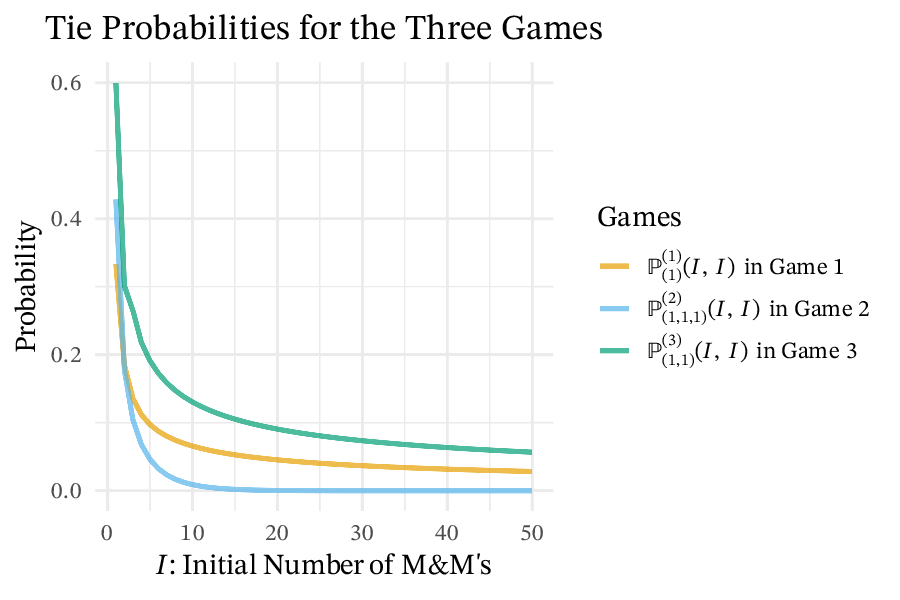}
    \caption{Tie Probabilities for the Three Games.}
    \label{fig: tie prob for three games}
\end{figure}

As an illustration, we computed tie probabilities in the case when each player starts with $I$ M\&M's and eats exactly one M\&M whenever a head is flipped; that is, $d=d_1=d_2=1$. For each game, we report both the exact value and its decimal approximation. These results are displayed in Table \ref{table:tie_prob_values} and Figure \ref{fig: tie prob for three games}. 

\begin{table}[!htbp]
\centering
\begin{tabular}{c cc cc cc}
\toprule
& \multicolumn{2}{c}{$\mathbb{P}^{(1)}_{(1)}(I_1, I_2)$ in Game 1} & \multicolumn{2}{c}{$\mathbb{P}^{(2)}_{(1,1,1)}(I_1, I_2)$ in Game 2} & \multicolumn{2}{c}{$\mathbb{P}^{(3)}_{(1,1)}(I_1, I_2)$ in Game 3} \\
\cmidrule(lr){2-3} \cmidrule(lr){4-5} \cmidrule(lr){6-7}
$(I_1,I_2)$
& Exact & Decimal
& Exact & Decimal
& Exact & Decimal \\
\midrule
$(1,1)$ & $1/3$ & $0.3333$ & $3/7$ & $0.4285$ & $3/5$ & $0.6000$ \\
$(2,2)$ & $5/27$ & $0.1851$ & $61/343$ & $0.1778$ & $113/375$ & $0.3013$ \\
$(3,3)$ & $11/81$ & $0.1358$ & $1759/16807$ & $0.1046$ & $22361/84375$ & $0.2650$ \\
\bottomrule
\end{tabular}
\caption{Comparison of Tie Probabilities for the Three Games.}
\label{table:tie_prob_values}
\end{table}

\subsection{Derivation of Tie Probabilities in Game 1}\label{subsection: derivation of game 1}

We now use Game 1 to illustrate our technique. Suppose a player eats $d > 0$ M\&M's if her coin lands heads. Inserting the step multiset $\mathcal{S} = [0, d] \times [0, d] \setminus [(0, 0) ]$ into the recurrence relation \eqref{eqs: general recurrence structure}, we have
\begin{equation*}
F(m,n)=
\begin{cases}
\begin{aligned}[t]
\frac{1}{3}\Bigl[F(m-1,n)&+F(m,n-1)\\
& +F(m-1,n-1)\Bigr] 
\end{aligned}
& \text{if } m\ge 1 \text{ and } n\geq 1  \\
0
& \text{if exactly one of } m,n \text{ is } \le 0  \\[6pt]

1
& \text{if } m\le 0 \text{ and } n\le 0.
\end{cases}
\end{equation*}

Define the generating function $A(x,y) = \sum_{m, n \geq 0}F(m, n)x^my^n$.  By Lemmas \ref{lemma: A(x,x) = 1 + sum F(m, n)x^m y^n} and \ref{lemma: rewrite F(m -a, n - b)}, we have
\begin{align*}
A(x,y) &\ = \ 1 +  \sum_{m \geq 1} \sum_{n \geq 1} F(m, n) x^m y^n \nonumber \\
&\ = \  1 + \frac{1}{3}\biggl[ \sum_{m \geq 1}\sum_{n \geq 1}F(m - d, n)x^m y^n   \nonumber \\
&\qquad \qquad + \sum_{m \geq 1} \sum_{n \geq 1}F(m, n - d)x^m y^n + \sum_{m \geq 1} \sum_{n \geq 1}F(m - d, n- d)x^m y^n \biggr]  \nonumber  \\
&\ = \  1 + \frac{1}{3}\Bigl(
    x^d\bigl(A(x,y)-1\bigr)
  + y^d\bigl(A(x,y)-1\bigr) \nonumber\\
&\qquad\qquad
  + \sum_{M=1}^{d}\sum_{N=1}^{d} x^M y^N
  + x^d y^d \bigl(A(x,y)-1\bigr)
\Bigr).
\end{align*}
For notational convenience, let $\alpha = x^d$ and $\beta = y^d$. After simple algebra, we find
\begin{align}
A(x,y) &\ = \  1 + \frac{\frac{1}{3}\sum^{d}_{M = 1} \sum^{d}_{N = 1} x^M y^N }{1 - \frac{1}{3}\left(\alpha + \beta + \alpha \beta \right)} \nonumber \\
&\ = \   1 +  \frac{1}{3}\sum_{M= 1}^{d} \sum_{N = 1}^{d}x^M y^N B\left(\alpha, \beta, \frac{1}{3} \right),
\label{eqs: A(x,y) in terms of B(alpha, beta, 1/3)}
\end{align}
where we recall that 
\begin{equation*}
B\left(\alpha, \beta, \frac{1}{3}\right) \ = \  \frac{1}{1 - \frac{1}{3}\left(\alpha + \beta + \alpha \beta\right)}.
\end{equation*}

We first study the expansion of $B\left(\alpha, \beta, \frac{1}{3}\right)$. After a simple application of the geometric series formula and multinomial expansion, we see that
\begin{equation*}
B\left(\alpha, \beta, \frac{1}{3}\right)  \ = \  \sum^{\infty}_{k = 0} \left(\frac{1}{3}\right)^{k} \sum_{\substack{i, j, \ell \geq 0 \\ i + j + \ell = k} } \frac{k!}{i!j! \ell!} \alpha^{i + \ell} \beta^{j + \ell}.
\end{equation*}
Let $m = i + \ell$ and $n  = j + \ell$. We then have
\begin{align*}
B\left(\alpha, \beta, \frac{1}{3}\right)  &\ = \  \sum_{m, n \geq 0} \alpha^m \beta^n \sum^{\min(m, n)}_{\ell = 0} \binom{m + n - \ell}{m - \ell, n - \ell, \ell}\left(\frac{1}{3}\right)^{m + n - \ell }\nonumber\\
&\ = \ \sum_{m, n \geq 0}\alpha^m \beta^n D\left(m, n, \frac{1}{3} \right).
\end{align*}
Now, we reparameterize $B\left(\alpha, \beta, \frac{1}{3}\right)$ by substituting $\alpha=x^d$ and $\beta=y^d$, and we see that
\begin{align}
B\left(\alpha, \beta, \frac{1}{3}\right) &\ = \ \sum_{m, n \geq 0} x^{dm} y^{dn} D\left(m, n, \frac{1}{3} \right). 
\label{eqs: expansion of B(x, y, 1/3)}
\end{align}

In Equation \eqref{eqs: expansion of B(x, y, 1/3)}, we have already computed the coefficients of $B\left(\alpha, \beta, \frac{1}{3}\right)$ in terms of $\alpha = x^{d}$ and $\beta = y^{d}$. To obtain the coefficients of $A(x, y)$, we insert  Equation \eqref{eqs: expansion of B(x, y, 1/3)} into Equation \eqref{eqs: A(x,y) in terms of B(alpha, beta, 1/3)}, and we have 
\begin{align*}
\label{eqs: generating function for Game 1 in expanded form}
A(x, y) &\ = \ 1 + \frac{1}{3}
\sum_{M=1}^d \sum_{N=1}^d x^M y^N\left( \sum_{m, n \geq 0} x^{dm}y^{dn} D\left(m, n, \frac{1}{3} \right)\right) \nonumber \\ 
&\ = \ 1 +  \frac{1}{3}
\sum_{M=1}^d \sum_{N=1}^d   \sum_{m, n \geq 0} x^{M + dm}y^{N + dn} D\left(m, n, \frac{1}{3} \right). 
\end{align*}
Let  $I_1 = M + dm$ and $I_2 = N + dn$ for $1\leq M \leq d$ and $1 \leq N \leq d$. We then have $I_1 - 1 = dm + (M - 1)$ and $I_2 -1 = dn + (N - 1)$, and $m, n$ exist by the division algorithm. In particular, we find that $m = \lfloor (I_1 - 1)/d \rfloor$ and $n = \lfloor (I_2 - 1) / d \rfloor$. Thus, we arrive at the following generating function
\begin{equation*} 
A(x, y) = 1 + \frac{1}{3}  \sum_{I_1, I_2 \geq 1} x^{I_1} y^{I_2} D\left(\left\lfloor \frac{I_1 - 1}{d} \right\rfloor, \left\lfloor \frac{I_2 - 1}{d} \right\rfloor, \frac{1}{3} \right). 
\end{equation*}
Indeed, we then find the coefficients of $A(x,y)$ as 
\[ 
\left[x^{I_1}y^{I_2}\right]A(x, y) \ = \ \frac{1}{3}D\left(\left\lfloor \frac{I_1 - 1}{d} \right\rfloor, \left\lfloor \frac{I_2 - 1}{d} \right\rfloor, \frac{1}{3} \right). 
\]


\subsection{Numerical Observations}
\label{subsec:simulations for three games}  
We first performed extensive computations to corroborate the correctness of our formulas. For each game, we computed the exact tie probability using both the recurrence relation \eqref{eqs: general recurrence structure} and the explicit formulas stated in Section \ref{subsection: results}. For Game 1, we computed $\mathbb{P}_{(d)}^{(1)}(I_1,I_2)$ for all $(I_1,I_2,d) \in [100]\times[100]\times[20]$, where $[N]=\{1,\dots,N\}$ for some positive integer $N$. For Game 2, we computed $\mathbb{P}_{(d,d_1,d_2)}^{(2)}(I_1,I_2)$ for all $(I_1,I_2,d,d_1,d_2) \in [100]\times[100]\times[20]\times[20]\times[20]$. For Game 3, we computed $\mathbb{P}_{(d_1,d_2)}^{(3)}(I_1,I_2)$ for all $(I_1,I_2,d_1,d_2) \in [100]\times[100]\times[20]\times[20]$. In every case, the values obtained from the two approaches were identical. Since the full output is too large to include in the paper, we provide a link to the code so that the reader can generate whatever comparisons they wish; see Appendix \ref{section: codes} for the code.

We now discuss two numerical observations. First, when both players begin with $I$ initial M\&M's, the number of coins each player tosses per round affects the behavior of the tie probability. In the original M\&M Game studied in \cite{badinski2017m}, as well as in our result for Game 1, the tie probability decreases with $I$. This decreasing behavior, however, does not hold uniformly when players are allowed to toss multiple coins.

\begin{figure}[!htbp]
    \centering    \includegraphics[width=0.9\textwidth]{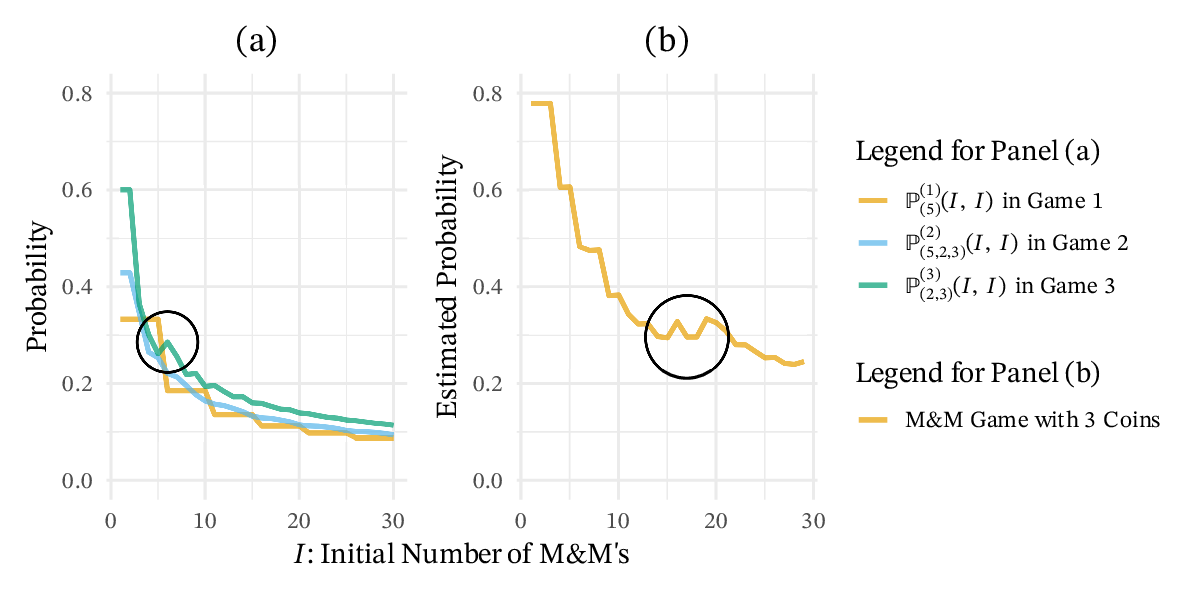}
    \caption{Effect of the Number of Coins on Tie Probabilities.}
    \label{fig: Starting number}
\end{figure}

In panel (a) of Figure \ref{fig: Starting number}, we compare $\mathbb{P}_{(5)}^{(1)}(I,I)$, $\mathbb{P}_{(5,2,3)}^{(2)}(I,I)$, and $\mathbb{P}_{(2,3)}^{(3)}(I,I)$ for $I \in [30]$. We make the following observation. While the tie probabilities in Games 1 and 2 are monotonically decreasing in $I$, the tie probability in Game 3 is no longer monotonically decreasing and exhibits small upward jumps at certain values of $I$. We circle one such region for illustration. For example, $\mathbb{P}_{(2,3)}^{(3)}(5,5)=0.262$, whereas $\mathbb{P}_{(2,3)}^{(3)}(6,6)=0.285$. In our simulations, we observed the same behavior for other choices of $(d_1,d_2)$. We also carried out Monte Carlo simulations for a three-coin version of the game and again found similar jumps. Panel (b) of Figure \ref{fig: Starting number} illustrates this, where we simulated the game in which a player eats $d_1=3$, $d_2=5$, and $d_3=10$ M\&M's when the corresponding coins land heads. These jumps suggest that multi-coin versions of the M\&M Game possess additional features that deserve further study.
 
\begin{figure}[!htbp]
    \centering
    \begin{minipage}{0.49\textwidth}
        \centering
        \includegraphics[width=\textwidth]{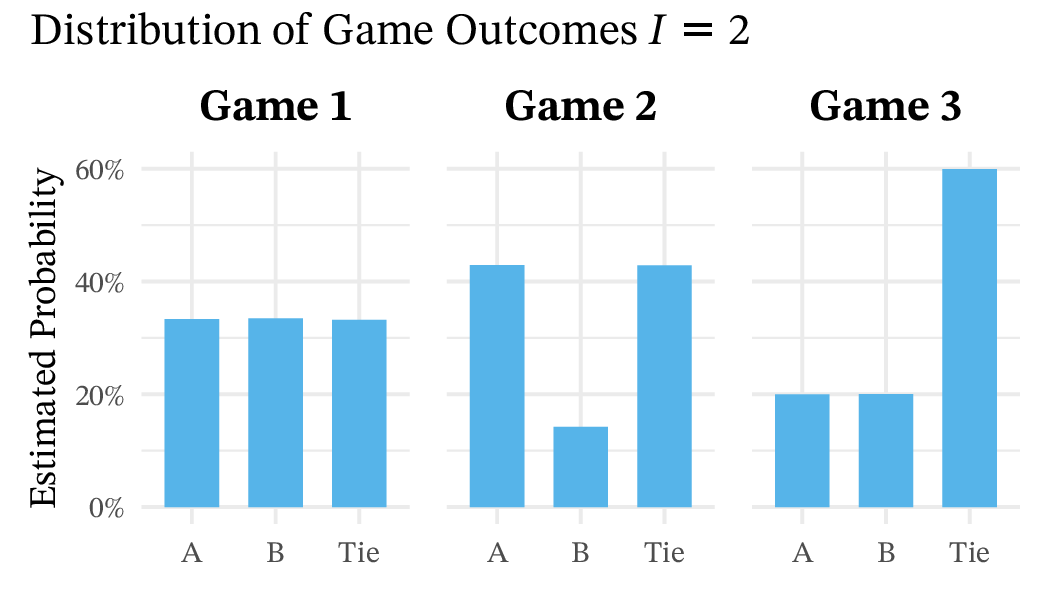}
    \end{minipage}
    \hfill
    \begin{minipage}{0.49\textwidth}
        \centering
        \includegraphics[width=\textwidth]{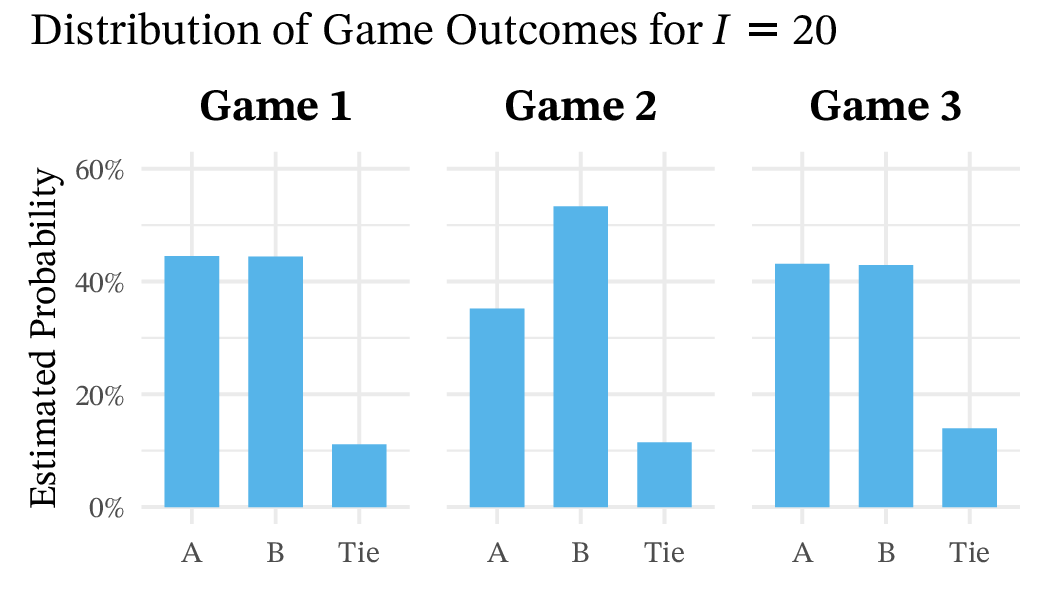}
    \end{minipage}
    \caption{Probability distributions of game outcomes for the three games when $I = 2$ and $I = 20$, based on Monte Carlo simulations. In both figures, \emph{A} denotes the probability that Player A wins the game, \emph{B} denotes the probability that Player B wins the game, and \emph{Tie} denotes the probability of a tie. }
    \label{fig:probability distributions of 3 games}
\end{figure}

Second, we found that when two players use Method 1 and Method 2, respectively, the probability distributions of the game outcomes change. In Figure~\ref{fig:probability distributions of 3 games}, we present Monte Carlo simulations for the three games when the starting number of M\&M’s is $I=2$ and $I=20$, and we plot the corresponding outcome distributions. On the horizontal axis, \emph{A} denotes the probability that Player A wins, \emph{B} denotes the probability that Player B wins, and \emph{Tie} denotes the probability of a tie. In these simulations, Game 1 uses $d=5$ for both players, Game 2 uses $d=5$ for Player A and $(d_1,d_2)=(2,3)$ for Player B, and Game 3 uses $(d_1,d_2)=(2,3)$ for both players. We find that the winning probabilities for Player A and Player B are symmetric in Games 1 and 3, but become asymmetric in Game 2. This asymmetry reinforces the phenomenon discussed earlier in Section~\ref{subsection: generalizing the MM Game}: even when two methods have the same expected number of M\&M’s removed per turn, they can still produce different probabilities of winning or losing, depending on which method each player uses.


\section{The M\&M Game Under Different Probabilities}
\label{sec: games under different probabilities}
\subsection{Extensions to Biased Coins}
The work in \cite{badinski2017m} and the extensions above focused only on the game using fair coins. It is natural to ask what happens when the two players use biased coins. We focus on the M\&M Game proposed by \cite{badinski2017m} with the modification that now each player uses a biased coin.  

By generalizing Theorem \ref{thm: probability tie finite sum}, we  easily have the following result. 
\begin{theorem}
\label{theorem: biased coins}
Consider two players who each start with $I$ M\&M's and flip one coin per round. On each round, a player eats one M\&M if her coin lands heads and eats none if it lands tails. Suppose each coin is biased, landing heads with probability $p$ and tails with probability $1-p$, where $0 < p < 1$. Then the probability of a tie is given by 
\[ 
\mathbb{P}(I, I) \ = \ \sum^{I - 1}_{k = 0} \binom{2I - k - 2}{k} \binom{2I - 2k - 2}{I - k - 1} \left(\frac{p}{2 - p} \right)^{k + 1}  \left(\frac{1 - p}{2 - p} \right)^{2I - 2k - 2}.  
\]
If $p = 0$, then $\mathbb{P}(I, I) = 0$. If $p = 1$, then $\mathbb{P}(I, I) = 1$. 
\end{theorem}
\begin{proof}
The cases for $p = 0$ and $p = 1$ are trivial. For $0 < p < 1$, the argument is similar to the proof given in \cite{badinski2017m}, except that we must recalculate the probability of each move after removing the null move $(0,0)$. Let $A$ denote the event that Player A flips heads and eats an M\&M, and define $B$ similarly for Player B. Then $A^c$ denotes the event that Player A does not eat an M\&M. Since removing the null move amounts to conditioning on the event $A\cup B$, we compute the probability that only one player eats an M\&M. The probability that only Player B eats an M\&M is given by 
\begin{align*} 
\mathbb{P}(A^c \cap B \mid A \cup B) 
&= \frac{\mathbb{P}(A^c \cap B)}{\mathbb{P}(A \cup B)} \\ 
&= \frac{p(1-p)}{p(1-p)+p(1-p)+p^2} \\
&= \frac{1-p}{2-p}.
\end{align*}
Because the two players' coins have the same probability $p$ of landing heads, we see $\mathbb{P}(A^c \cap B \mid A \cup B)  = \mathbb{P}(A\cap B^{c}\mid A \cup B)$. The probability that both players eat an M\&M can be computed by a similar argument and is given by 
\[
\mathbb{P}(A\cap B \mid A\cup B)
=
\frac{p^2}{2p(1-p)+p^2}
=
\frac{p}{2-p}.
\]
The rest of the combinatorial argument is the same as in \cite{badinski2017m}.
\end{proof}

\begin{figure}[h]
    \centering
    \includegraphics[width=0.8\linewidth]{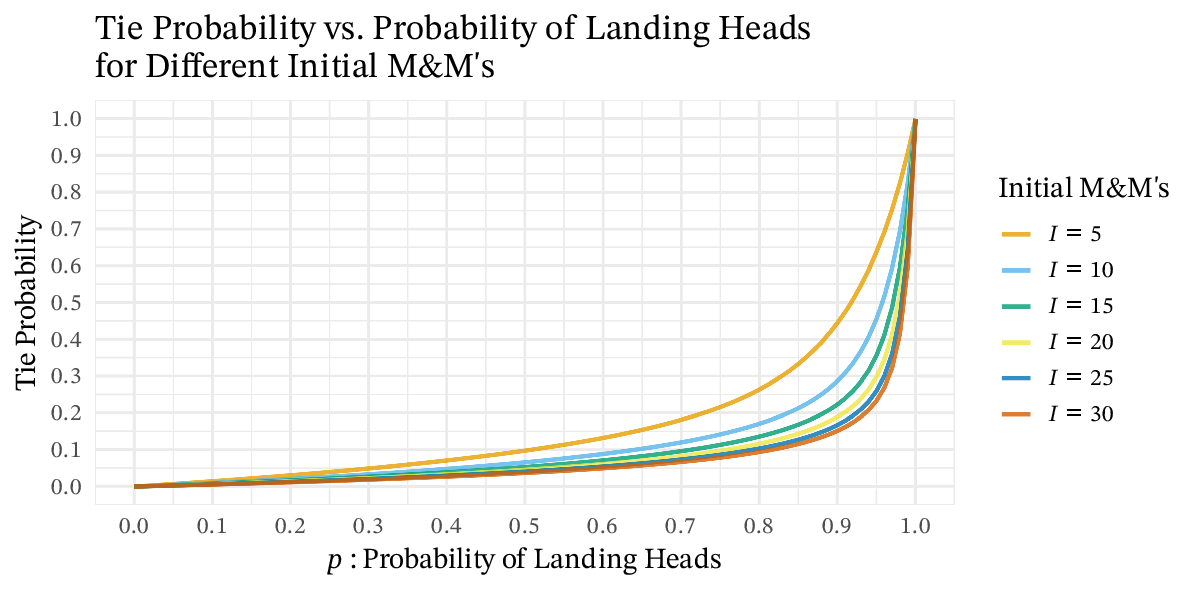}
    \caption{Tie Probability as a Function of $p$  for Different Initial Numbers of M\&M's.}
    \label{fig: tie function}
\end{figure}

Using Theorem \ref{theorem: biased coins}, we computed the tie probabilities for the game when each player starts with $I$ M\&M's and flips a biased coin landing heads with probability $p$ from the set $\{0, 0.1, \dots, 1\}$. The results are displayed in Figure \ref{fig: tie function}. As the figure shows, the tie probability is an increasing function of $p$, with the increase becoming particularly steep when $p \geq 0.9$. This behavior is intuitive: as $p$ increases, both players are more likely to eat an M\&M in each round, which raises the chance that they finish their last M\&M in the same round.  

\subsection{The M\&M Game with Evolving Coin Probabilities} 
We propose a state-dependent version of the M\&M Game in which the probability that a coin lands heads evolves according to an exponential model. In this version of the game, the probability that a player’s coin lands heads depends on the number of M\&M's that player has remaining: in particular, the fewer M\&M's a player has left, the greater the probability that the player’s coin lands heads. This reflects the idea that the closer a player is to depletion (that is, death), the more likely the player is to consume another M\&M. Because the probability of landing heads changes over the course of the game, we call this version of the M\&M Game the \emph{Evolving Coin Probabilities Game}.

Modifying the setup of the M\&M Game in \cite{badinski2017m}, we define the Evolving Coin Probabilities Game  as follows. Two players each begin with $I$ M\&M's, and in each round each player independently flips a coin. A player consumes one M\&M if the coin lands heads and consumes none if it lands tails. The probability that a player’s coin lands heads in a given round is
\begin{equation}
    \mathbb{P}(\text{head})\ =\ 1 - e^{-\lambda(I - q + 1)},\label{eqs: exponential coin clip}
\end{equation}
where $\lambda > 0$ is the rate parameter controlling the rate of increase, and $q$ is the number of M\&M's currently remaining for that player. The rules determining when a player wins and when the game ends in a tie are the same as the rules given in Section \ref{subsection: background of the MM Game}.

Given the initial number of M\&M's $I$ and the rate parameter $\lambda$, we see from Equation \eqref{eqs: exponential coin clip} that the probability of landing heads depends on the difference $I-q$, which is the number of M\&M's the player has already consumed. As a player consumes more M\&M's, the probability in Equation \eqref{eqs: exponential coin clip} increases toward $1$. We study the tie probability for the Evolving Coin Probabilities Game. Let $\mathbb{P}(\text{tie}; \lambda, I)$ denote the probability of a tie in this game with rate parameter $\lambda$ and initial number of M\&M's $I$. In what follows, we conduct extensive simulations to study $\mathbb{P}(\text{tie}; \lambda, I)$.

\subsubsection{Simulation Results}
\label{subsec: simulation results for evolving probability}
\begin{figure}[!htbp]
    \centering
    \begin{minipage}[t]{0.495\linewidth}
        \centering
        \includegraphics[width=\linewidth]{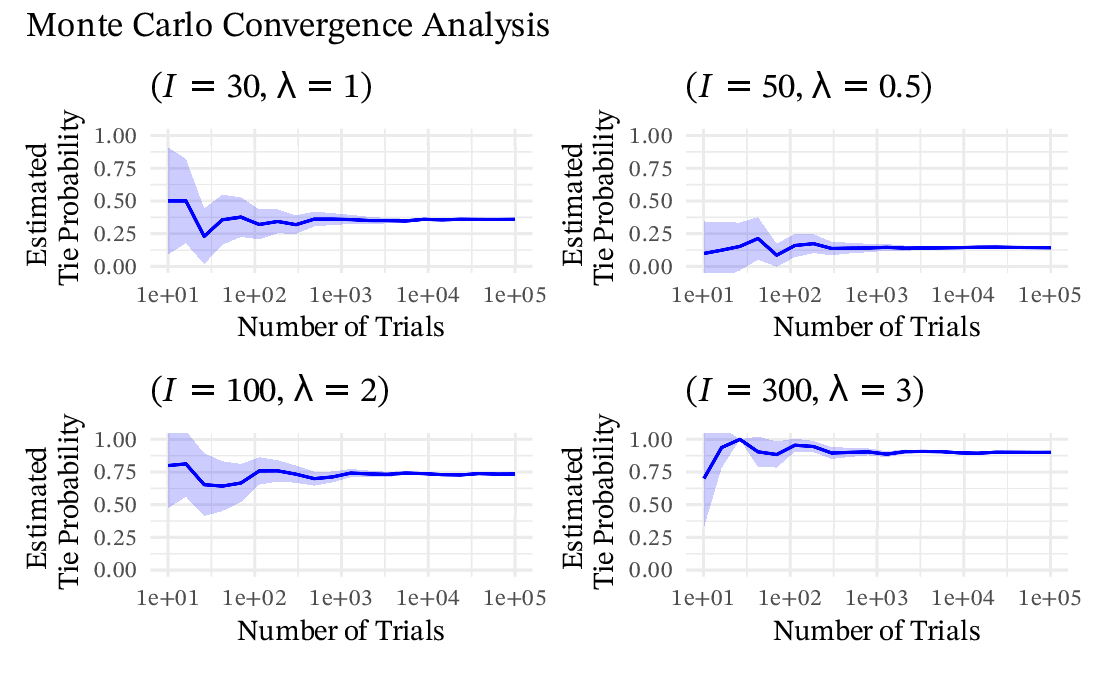}
        \caption{Monte Carlo Convergence Analysis for Selected Values of $(\lambda, I)$.}
        \label{fig:monte_carlo_convergence} 
    \end{minipage}\hfill
    \begin{minipage}[t]{0.495\linewidth}
        \centering
        \includegraphics[width=\linewidth]{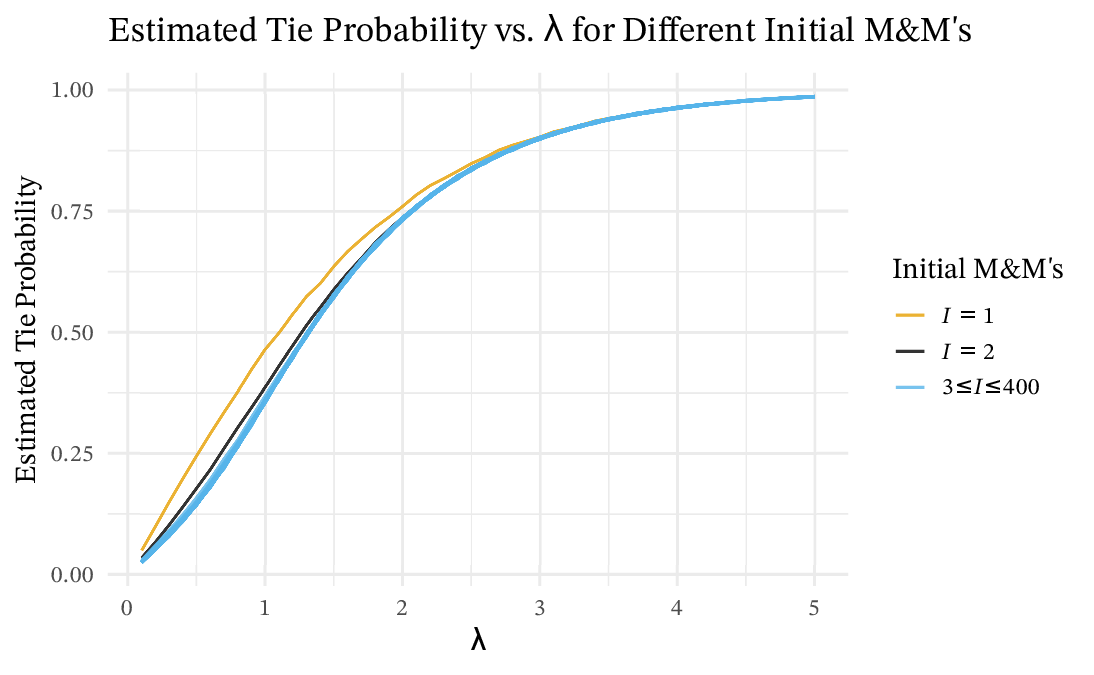}
        \caption{Estimated Tie Probabilities as a Function of the Rate Parameter  $\lambda$.}
        \label{fig:probabilityExp_lambda}
    \end{minipage}
\end{figure}

We performed extensive Monte Carlo simulations to examine the influence of parameters $\lambda$ and $I$ on $\mathbb{P}(\text{tie}; \lambda, I)$. Specifically, we simulated the game for all parameter pairs $(\lambda,I)$ with $\lambda \in \{0.1,0.2,\dots,5\}$ and $I \in \{1,2,\dots,400\}$, using step size $0.1$ in $\lambda$. For each choice of $(\lambda, I)$, we performed 100,000 simulation trials. The results are displayed in Figure \ref{fig:probabilityExp_lambda}. In the figure, the yellow curve corresponds to the case $I=1$, the black curve corresponds to $I=2$, and the collection of blue curves corresponds to $I \in \{3,\dots,400\}$. To facilitate the discussion, we use $\mathbb{P}^{(\mathrm{mc})}_{T}(\text{tie}; \lambda, I)$ to denote the Monte Carlo estimate of the tie probability for a fixed rate parameter $\lambda$ and an initial number of M\&M's $I$, and the subscript $T$ is used to denote that the Monte Carlo estimate is based on $T$ simulation trials. 

To evaluate the quality of our Monte Carlo simulations, we computed the $99\%$ confidence intervals of $\mathbb{P}^{(\mathrm{mc})}_{T}(\text{tie}; \lambda, I)$ for different numbers of trials $T$ over the range $10\leq T\leq 100{,}000$. The $99\%$ confidence interval is computed according to the  binomial proportion confidence interval
\begin{equation*} 
\mathbb{P}_{T}^{(\text{mc})}(\text{tie}; \lambda, I)\ \pm \  \frac{2.58}{\sqrt{T}} \sqrt{ \mathbb{P}_{T}^{(\text{mc})}(\text{tie}; \lambda, I)\left(1 - \mathbb{P}_{T}^{(\text{mc})} (\text{tie}; \lambda, I) \right)}, 
\end{equation*}
where $T$ denotes the number of simulation trials and $2.58$ is the $z$-score corresponding to a $99\%$ confidence level for the standard normal distribution. For a detailed overview on the construction of confidence intervals for Monte Carlo simulations, see \cite{mcbook}. Selected examples are shown in Figure \ref{fig:monte_carlo_convergence}, where the blue curve represents the Monte Carlo estimates of the tie probabilities and the shaded region represents the corresponding confidence intervals. As the number of simulation trials $T$ increases, the confidence intervals steadily narrow and become negligible once the number of simulation trials reaches 100,000. This indicates that the estimate $\mathbb{P}_{T}^{\text{(mc})}(\text{tie}; \lambda, I)$ converges sufficiently to the true tie probability $\mathbb{P}(\text{tie}; \lambda, I)$. 

Based on the simulation results, we make two key observations about the Evolving Coin Probabilities Game. First, the tie probability appears to be increasing in $\lambda$: as $\lambda$ grows, the probability of a tie also increases. This is because larger values of $\lambda$ cause the probability of flipping heads to rise more rapidly over the course of the game, thereby increasing the likelihood that both players consume their M\&M's in the same round and hence eat their last M\&M's in the same round. By contrast, smaller values of $\lambda$ lead to a lower probability of flipping heads, which reduces the likelihood that both players consume their M\&M's in the same round and thereby lowers the tie probability.

Second, the tie probabilities for different values of $I$ appear to align remarkably well. As illustrated in Figure \ref{fig:probabilityExp_lambda}, the curves for different initial numbers of M\&M's $I$ begin to cluster together once $I \geq 3$. This indicates that, for moderate to large values of $I$, the tie probability is nearly independent of $I$. This behavior stands in sharp contrast to the M\&M Game with biased coins (see Figure \ref{fig: tie function}), where for a fixed probability of landing heads $p$, the tie probability decreases as $I$ increases. In that setting, the initial number of M\&M's $I$ plays a more visible role in determining the game outcome. In the Evolving Coin Probabilities Game, however, the number of initial M\&M's $I$ exerts comparatively little influence on the tie probability once $I$ is moderately large. Instead, the rate parameter $\lambda$ emerges as the primary factor governing the tie probabilities.

We performed a Pearson correlation analysis \cite{rice2006mathematical} to support our observations. The results show a strong positive correlation between the rate parameter $\lambda$ and the tie probability, with a correlation coefficient of $0.9251$. In contrast, the correlation between the initial number of M\&M's $I$ and the tie probability is negligible, with a correlation coefficient of $-0.0006$. Since the correlation between the tie probability and $\lambda$ is close to $1$, while the correlation between the tie probability and $I$ is close to $0$, these findings provide further empirical support for our observation that $\lambda$ is the dominant factor governing the tie probability, whereas the initial number of M\&M's has minimal influence on the outcome of the game.


\subsubsection{Modeling Tie Probabilities with a Gompertz Curve}\label{subsec:gompertz_fit}
In Figure~\ref{fig:probabilityExp_lambda}, we observed that the simulated tie probabilities $\mathbb{P}_{T}^{(\mathrm{mc})}(\textup{tie};\lambda,I)$ nearly collapse onto a single curve for all $I \geq 3$, which suggests that the influence of $I$ on the tie probability becomes negligible for moderately large $I$. The correlation analysis similarly supports the conclusion that the tie probability is nearly independent of $I$. Based on these observations, it is natural to ask whether a probability model  depending only on $\lambda$, denoted $G(\lambda)$, could approximate the tie probability $\mathbb{P}(\text{tie}; \lambda, I)$ well for moderately large values of $I$. Since the simulated probability curves exhibit slow initial growth, followed by a phase of more rapid increase, and then level off toward an asymptote, a Gompertz curve appears to be an appropriate choice for modeling this behavior.

We now present a brief overview of the Gompertz curve. The Gompertz curve is commonly used to model growth processes characterized by an initial slow phase, followed by a phase of rapid increase, and eventual saturation; see  \cite{dhar2018comparison} for a detailed discussion of Gompertz curves. Let $G(\lambda)$ denote the tie probability for the Evolving Coin Probabilities Game as predicted by the Gompertz model. The Gompertz model is given by
\begin{equation}
\label{eqs: general form of gompertz model}
     G(\lambda) \ = \ L   \exp\left( -h   \exp(-\lambda_0 \lambda) \right). 
\end{equation}
In Equation \eqref{eqs: general form of gompertz model}, the parameter $L$ is the upper asymptote, representing the limiting value of the growth process.  Clearly, we see that $G(\lambda) \to L$ as $\lambda \to \infty$.  The parameter $h$  determines the horizontal placement of the curve. Finally, $\lambda_0$ is the growth rate, determining how quickly the curve rises toward its asymptote.  

We now seek to fit a Gompertz model to the simulated tie probabilities. Because the parameters $(L,h,\lambda_0)$ are unknown, we estimate them by minimizing the sum of squared errors between the simulated probabilities and the Gompertz model given in Equation \eqref{eqs: general form of gompertz model}. Since the resulting optimization problem is nonlinear, a closed-form solution for the optimal parameters is generally unavailable. Accordingly, we employ a numerical least-squares procedure, specifically the Levenberg--Marquardt algorithm, to obtain estimates of $(L,h,\lambda_0)$ \cite{marquardt1963algorithm}. The fitting procedure is summarized in Algorithm~\ref{algorithm: gompertz curve fit}.

\begin{algorithm}[!htpb]
\caption{Parameter Estimation Procedure for Gompertz Model}
\label{algorithm: gompertz curve fit}
\begin{algorithmic}[1]
\STATE \textbf{Input:}   $\Lambda = \{\lambda_i\}_{i = 1}^{50}$ for $\lambda_i = i/10$, number of Monte Carlo trials $T = 100{,}000$, and initial numbers of M\&M's $I \in \{3,\dots,400\}$.\STATE \textbf{Output:} 
Estimated parameters for the Gompertz model defined in Equation \eqref{eqs: general form of gompertz model}.
\STATE Create an empty list $\chi$.  
\FOR{each $\lambda_i \text{ in } \Lambda$}
\STATE Simulate $\mathbb{P}_{T}^{(\text{mc})}(\text{tie}; \lambda_i, I)$ for all values of $I \in \{3, \dots, 400\}$.
\STATE Average the values of $\mathbb{P}_{T}^{(\text{mc})}(\text{tie}; \lambda_i, I)$ over $I$ to obtain 
    \[\overline{\mathbb{P}^{\text{(mc)}}_{T}}(\lambda_i)
    \ = \ 
    \frac{1}{398}\sum_{I=3}^{400}
    \mathbb{P}^{(\text{mc})}_{T}(\textup{tie};\lambda_i, I).
    \]
    \STATE Append the value $\left(\lambda_i, \overline{\mathbb{P}^{\text{(mc)}}_{T}}(\lambda_i)\right)$ to  $\chi$.
\ENDFOR

\STATE Use the Levenberg--Marquardt algorithm to solve the least-squares problem     
\[
\min_{L,h,\lambda_0}
\sum_{\left(\lambda_i, \overline{\mathbb{P}^{\text{(mc)}}_{T}}(\lambda_i) \right) \in \chi} 
\left(
\overline{\mathbb{P}^{\text{(mc)}}_{T}}(\lambda_i) - L\exp\left(-h \exp\left(-\lambda_0 \lambda_i\right) \right)
\right)^2, 
\]
and obtain the estimated parameters $(\widehat{L},\widehat{h},\widehat{\lambda_0})$ that  minimize the quantity above.

\STATE \textbf{Return:} $(\widehat{L}, \widehat{h}, \widehat{\lambda_0})$.  
\end{algorithmic}
\end{algorithm}

\begin{figure}[!htpb]
    \centering
    \begin{minipage}{0.49 \linewidth}
        \centering
        \includegraphics[width=\linewidth]{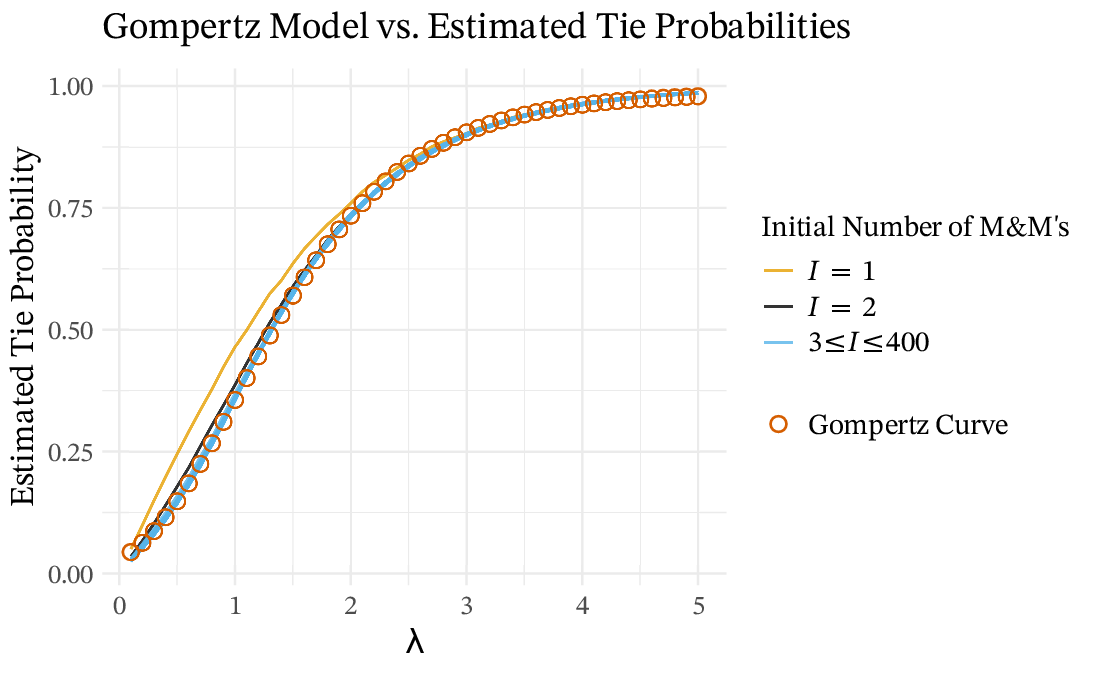}
        \caption{Modeling Tie Probabilities \\  via Gompertz Model.}
        \label{fig:Gompertz Model} 
    \end{minipage} 
    \begin{minipage}{0.49 \linewidth}
        \centering
        \includegraphics[width=\linewidth]{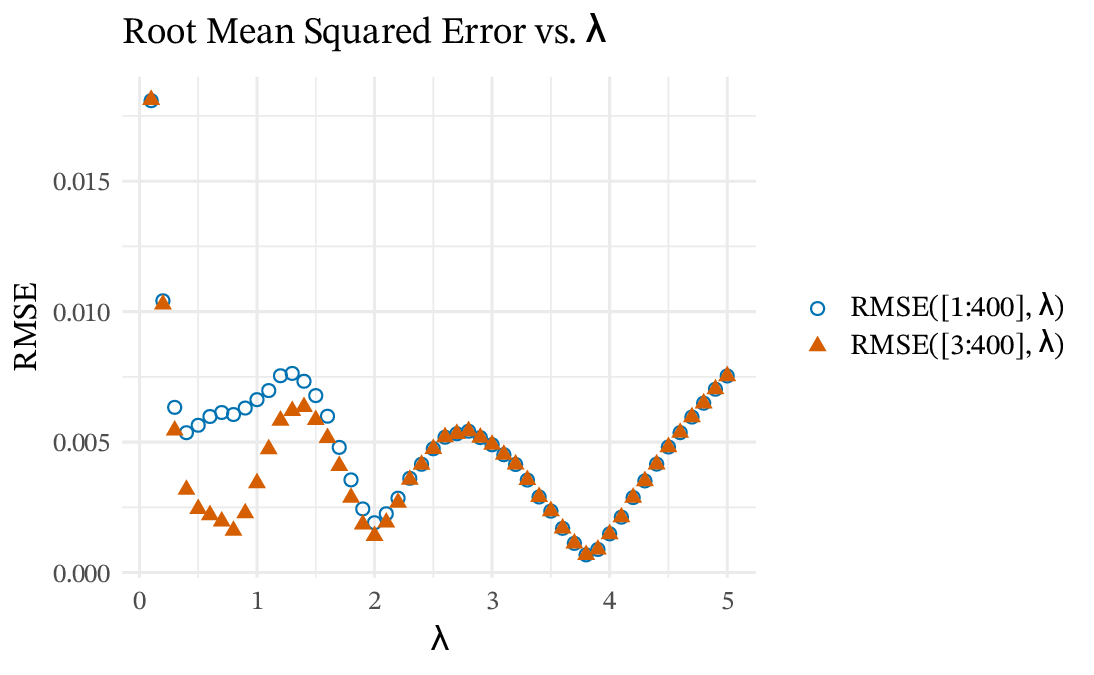}
        \caption{RMSE of the Fitted Gompertz Model.}
        \label{fig:RMSE} 
    \end{minipage}
\end{figure}

By applying Algorithm \ref{algorithm: gompertz curve fit}, we obtain the estimated parameters $\widehat{L} \approx 0.986, \widehat{h} \approx 3.525$, and $ \widehat{\lambda_0} \approx 1.241$. Our final model is then given by
\begin{equation} 
\label{eqs: fitted gompertz model}
\widehat{G}(\lambda) \ \approx \ 0.986 \exp\!\left(-3.525 \exp(-1.241 \lambda)\right). 
\end{equation}
Our model aligns with our expectations. In particular, the estimated parameter $\widehat{L}\approx0.986$ is very close to 1, which suggests that the fitted Gompertz curve $\widehat{G}(\lambda)$ has a horizontal asymptote near 1. This agrees with our numerical observation that the tie probabilities approach 1 as $\lambda$ becomes sufficiently large. In Figure \ref{fig:Gompertz Model}, we plot the fitted Gompertz model $\widehat{G}(\lambda)$ together with the simulated curves shown earlier in Figure~\ref{fig:probabilityExp_lambda}. The tie probabilities predicted by Equation \eqref{eqs: fitted gompertz model}, represented by red circles in Figure~\ref{fig:Gompertz Model}, capture the overall trend of the tie probabilities and agree well with $\mathbb{P}_{T}^{\text{(mc)}}(\text{tie};\lambda,I)$ for all $I \geq 3$. As expected, the Gompertz model does not fully predict the cases $I=1$ and $I=2$, since for such small initial values the tie probability still depends on $I$.

For each value of $\lambda$, we also assessed the goodness of fit of the Gompertz model to the simulated probability curves $\mathbb{P}_{T}^{(\text{mc})}(\text{tie}; \lambda, I)$. We used the root mean squared error (RMSE) as our measure of fit. One advantage of RMSE is that it is reported on the same scale as the original quantity being modeled. In our setting, for example, an RMSE of $0.001$ indicates that the fitted model typically differs from the simulated tie probability by about $0.001$. Specifically, we computed two RMSE values: $\text{RMSE}([1:400], \lambda)$, based on all simulated tie probabilities from Section~\ref{subsec: simulation results for evolving probability}, and $\text{RMSE}([3:400], \lambda)$, based only on the probabilities for $3 \leq I \leq 400$, where the simulated curves effectively collapse onto a single curve. These quantities are computed as follows:
\begin{equation*} 
\text{RMSE}([1:400], \lambda)  \ = \   \sqrt{ \frac{1}{400}\sum_{I \in [1:400]} \left( \mathbb{P}_{T}^{\text{(mc)}}(\text{tie};\lambda, I) - \widehat{G}(\lambda) \right)^2 }  ,
\end{equation*}
and
\begin{equation*} 
\text{RMSE}([3:400], \lambda) \ = \    \sqrt{ \frac{1}{398}\sum_{I \in [3:400]} \left( \mathbb{P}_{T}^{\text{(mc)}}(\text{tie};\lambda, I) - \widehat{G}(\lambda)\right)^2 }. 
\end{equation*}
We plot these quantities in Figure \ref{fig:RMSE}. We observe that $\text{RMSE}([1:400], \lambda)$ is larger than $\text{RMSE}([3:400], \lambda)$ for $\lambda \in [0.5,2]$. This is consistent with our earlier observation in Figure \ref{fig:Gompertz Model} that the probability curves for $I=1$ and $I=2$ do not align closely with the Gompertz model in this interval. Outside this interval, however, the two RMSE curves agree closely. We also observe that the RMSE is largest for $\lambda < 0.5$, suggesting a greater discrepancy between our fitted model and the tie probability in this region. One possible explanation is that, for small values of $\lambda$, the growth rate predicted by the Gompertz model does not match the true rate at which the tie probability changes with $\lambda$. Nevertheless, the Gompertz model still approximates the tie probability accurately to two decimal places, indicating that it remains a useful model for approximating tie probabilities.


While the Gompertz model captures the overall trend in the tie probabilities reasonably well, it still has several limitations. First, our fitted model only accounts for tie probabilities when $I \geq 3$. Second, the Gompertz curve imposes a particular asymmetric form: slow initial growth, followed by more rapid increase, and finally gradual saturation. Although our simulated tie probabilities are broadly consistent with this pattern, the rate of growth is not equally well captured across all values of $\lambda$. This likely contributes to the wave-shaped behavior of the RMSE curve. For these reasons, it may be worthwhile to explore alternative models that better capture these features, and we leave this for future exploration.

\section{Conclusion and Future Work}\label{sec: conclusion}
We investigated several extensions of the M\&M Game, focusing on how modified game rules and probability distributions of coin flips affect the tie probability. We now discuss future work for the M\&M Game. 

In Section \ref{sec: Comparison of Three Games}, we significantly generalized the M\&M Game, though many further generalizations remain possible. First, it would be worthwhile to study the M\&M Game in the more general setting of $P$ players rather than only two players. Second, we studied a generalized version of the M\&M Game in which a player tosses two coins in each round, consuming $d_1$ M\&M's if the first coin lands heads and $d_2$ M\&M's if the second coin lands heads. It is interesting to study the tie probability when a player can flip $r$ coins in a given round. In particular, what is the tie probability when, on each turn, a player may consume $d_1, d_2, \dots, d_r$ M\&M's, where the probability of the $i^{\text{th}}$ coin landing heads is $p_i$? We leave these more general cases for future investigation.

Besides studying the tie probability for generalized versions of the M\&M Game, another unexplored question in both \cite{badinski2017m} and this paper concerns the length of the game: what is the expected number of turns for the game to terminate? During the course of writing this paper, we examined this problem, and we believe it could similarly be solved using a recurrence relation. We briefly sketch the idea. Let $E(m,n)$ denote the expected number of turns for the game to terminate when Player A has $m$ M\&M's and Player B has $n$ M\&M's. Previously, we defined $\mathcal{S} = (\mathcal{A} \times \mathcal{B}) \setminus [(0, 0)]$, omitting $(0,0)$ because we were concerned only with game turns in which the numbers of M\&M's possessed by the two players change. When computing the expected number of turns, however, we must also account for the possibility that neither player eats an M\&M in a given round. Therefore, we need to consider $\mathcal{S} \cup [(0, 0)]$, the multiset of all possible moves in a round. Using $\mathcal{S} \cup [(0, 0)]$, we can obtain the following recurrence relation:  
\begin{align*} 
E(m, n) &\ = \ 1 +  \frac{1}{|\mathcal{S}| + 1}\left( \sum_{(a, b)\in \mathcal{S} \cup [(0, 0)] }E(m - a, n - b)\right) \nonumber \\ 
&\ = \  1 +  \frac{1}{|\mathcal{S}| + 1}\left( \sum_{(a, b) \in \mathcal{S}} E(m - a, n - b) + E(m, n)\right) \nonumber \\ 
&\ = \  \frac{|\mathcal{S}| + 1}{|\mathcal{S}|} + \frac{1}{|\mathcal{S}|} \sum_{(a, b) \in \mathcal{S}} E(m - a, n - b) .  
\end{align*}
We observe $E(0,0)=0$, since the expected number of turns for the game to terminate is $0$ when the game is already at $(0,0)$. We also note $E(m,0)=E(0,n)=0$ for $m,n>0$, since the game terminates immediately once one player has exhausted all of her M\&M's. The technique used in this paper, namely computing the generating function and then extracting coefficients, can similarly be applied to derive an exact formula for the expected number of turns in the M\&M Game. We leave this problem, along with its possible extensions, as a direction for future work.

In Section \ref{sec: games under different probabilities}, we studied the M\&M Game under modified coin-flip probability distributions. In our study, we have provided only empirical evidence that $\lambda$ is the dominant parameter governing the tie probability and that, to a large extent, the tie probability is nearly independent of $I$. It would be valuable to develop a rigorous theoretical explanation for this phenomenon, as well as to obtain a closed-form expression for the tie probability in the Evolving Coin Probabilities Game. In addition, it would be interesting to study the M\&M Game under alternative coin-flip distributions.

\newpage 
\appendix

\section{Code}
\label{section: codes}
All simulations in the paper were conducted in Python. For readers interested in replicating our studies, the code is available at \url{https://github.com/MatthewYilong/MMGame}.

\section{Proof of Lemmas}\label{appendix: Proof of Lemmas}
\subsection{Proof of Lemma \ref{lemma: A(x,x) = 1 + sum F(m, n)x^m y^n}} 
A simple computation and an application of the recurrence relation \eqref{eqs: general recurrence structure} yield
\begin{align*}
A(x, y)  &\ = \  \sum_{m \geq 0} \sum_{n \geq 0} F(m, n)x^m y^n \nonumber \\
&\ = \   F(0,0) + \sum_{n>0} F(0,n)y^n + \sum_{m>0} F(m,0)x^m + \sum_{m>0}\sum_{n>0} F(m,n)x^m y^n \nonumber \\
&\ = \  1 + 0 + 0 + \sum_{m>0}\sum_{n>0} F(m,n)x^m y^n \nonumber \\
  &\ = \ 1 + \sum_{m \geq 1}\sum_{n \geq 1} F(m,n)x^m y^n,
\end{align*}
where $\sum_{n>0} F(0,n)y^n$ and $\sum_{m > 0}F(m, 0)x^m$ vanish as a result of the initial conditions of recurrence relation \eqref{eqs: general recurrence structure}.
\qed 

\subsection{Proof of Lemma \ref{lemma: rewrite F(m -a, n - b)}}
We prove the first formula. Consider $\sum_{m, n \geq 1} F(m - a, n - b)x^m y^n$ for $a, b \geq 1$. If we let $i = m - a$ and $j = n - b$, then we have
\begin{align*}
\sum_{m, n \geq 1} F(m - a, n - b)x^m y^n &\ = \ \sum_{i + a \geq 1} \sum_{j + b \geq 1} F(i, j)x^{i + a} y^{j + b} \nonumber \\
&\ = \ x^{a}y^{b} \sum_{i \geq - (a - 1)}\sum_{j \geq - (b- 1)}F(i, j) x^i y^j \nonumber \\
&\ = \ x^{a}y^{b} \left(\sum_{i = - (a- 1)}^{0} \sum_{j = - (b - 1)}^{0}F(i,j)x^i y^j + \sum_{i, j \geq 1} F(i, j)x^i y^j \right) \nonumber \\
&\ = \  \sum^{0}_{i = -(a - 1)}\sum^{0}_{j = -(b - 1)} F(i,j)x^{i + a} y^{j + b} + x^{a}y^{b} \sum_{i, j \geq 1} F(i, j)x^i y^j \nonumber \\
&\ = \  \sum^{0}_{i = -(a - 1)}\sum^{0}_{j = -(b - 1)} x^{i + a} y^{j + b} + x^{a}y^{b} \sum_{i, j \geq 1} F(i, j)x^i y^j \nonumber \\
&\ = \ \sum_{M = 1}^{a} \sum_{N = 1}^{b}x^M y^N + x^{a}y^{b} \left( A(x, y) - 1\right).
\end{align*}
In the penultimate step, we see that $F(i, j) = 1$ if both $i \leq 0$ and $j \leq 0$, so we have that $\sum^{0}_{i = -(a - 1)}\sum^{0}_{j = -(b - 1)} F(i,j)x^{i + a} y^{j + b}$ equals $ \sum^{0}_{i = -(a - 1)}\sum^{0}_{j = -(b - 1)} x^{i + a} y^{j + b}$. In the last step, we let $M = i + a$ and $N = j + b$ for the term $ \sum^{0}_{i = -(a - 1)}\sum^{0}_{j = -(b - 1)} x^{i + a} y^{j + b}$, and we apply Lemma \ref{lemma: A(x,x) = 1 + sum F(m, n)x^m y^n} to the term $\sum_{i, j \geq 1}F(i, j)x^i y^j$.

We next prove the formula for $\sum_{m, n \geq 1}F(m - a, n)x^m y^n$. Let $i = m - a$. We then have
\begin{align*}
\sum_{m, n \geq 1} F(m - a, n)x^m y^n &\ = \ \sum_{i + a \geq 1} \sum_{n \geq 1} F(i, n)x^{i + a} y^n  \nonumber \\
&\ = \ x^a \sum_{i \geq - (a - 1)} \sum_{n \geq 1} F(i, n) x^i y^n \nonumber \\
&\ = \ x^a  \left\{\sum_{i = - (a - 1)}^{0} \sum_{n \geq 1}F(i, n)x^i y^n +  \sum_{i \geq 1} \sum_{n \geq 1} F(i, n)x^i y^n  \right\} \nonumber \\
&\ = \  x^{a}\left( A(x, y) - 1 \right),
\end{align*}
where the first sum vanishes due to our initial conditions in recurrence relation \eqref{eqs: general recurrence structure} and we apply Lemma \ref{lemma: A(x,x) = 1 + sum F(m, n)x^m y^n} to the second term. The proof of the formula for $\sum_{m, n \geq 1}F(m, n - b)x^m y^n$ is exactly the same, so we omit it here.
\qed


\section{Proofs of Results in Section \ref{subsection: results}}\label{appendix: Proof of Results for Game 2 and Game 3}

\subsection{Proof of Theorem \ref{thm: connection between Game 1 and original MM Game}} 
Observe that 
\begin{align*}
\mathbb{P}_{(d)}^{(1)}\left(dI, dI\right) &\ =  \ \frac{1}{3}D\left( \left\lfloor \frac{dI - 1}{d}\right\rfloor,  \left\lfloor \frac{dI - 1}{d}\right\rfloor, \frac{1}{3}  \right) \\ 
&\ = \ \frac{1}{3} D\left( \left\lfloor I - \frac{1}{d} \right\rfloor, \left\lfloor I - \frac{1}{d} \right\rfloor, \frac{1}{3}  \right) \\ 
&\ = \  \frac{1}{3}D\left(I + \left\lfloor - \frac{1}{d} \right \rfloor, I + \left\lfloor - \frac{1}{d} \right \rfloor, \frac{1}{3} \right) \\
&\ = \ \frac{1}{3}D\left(I - 1,  I - 1, \frac{1}{3}\right) \\ 
&\ = \ \mathbb{P}^{\mathrm{original}}(I, I).
\end{align*} 
In the penultimate step, we see that $\frac{1}{3} D(I-1,I-1,\frac{1}{3})$ is exactly the tie probability in the original M\&M Game. This proves the desired result.
\qed 

\subsection{Proof of Theorem \ref{thm: finite sum for Game 2}} 
We analyze Game 2. Let $d, d_1, d_2 \geq 1$, so we have the following recurrence relation
\begin{equation*}
F(m,n) =
\begin{cases}
\begin{aligned}
\frac{1}{7}\biggl[&F(m-d,n-d_1)
+ F(m-d,n-d_2) \\ 
& +  F(m-d,n-d_1-d_2) 
+ F(m-d,n) \\
&+  F(m,n-d_1)  
+ F(m,n-d_2) + F(m, n - d_1 - d_2)  \biggr]
\end{aligned}
& \text{ if } m,n \geq 1, \\
0 & \text{ if } m\leq0,n \geq 1  \\
& \text{ or } n\leq0,m \geq 1  \\
1 & \text{ if } m, n \leq 0.
\end{cases}
\end{equation*}

Define the generating function $A(x, y) = \sum_{m, n \geq 0} F(m, n) x^m y^n$.  We have the following
\begin{align*}
A(x, y) &\ = \  1  +  \frac{1}{7} \sum_{m \geq 1} \sum_{n \geq 1} F(m - d, n - d_1)x^{m}y^{n}  + \frac{1}{7} \sum_{m \geq 1} \sum_{n \geq 1} F(m -d, n-d_2)x^{m}y^{n}  \nonumber \\
&\qquad \qquad  +  \frac{1}{7} \sum_{m \geq 1} \sum_{n \geq 1}F(m - d, n - d_1 -d_2)x^{m}y^{n} +  \frac{1}{7} \sum_{m \geq 1} \sum_{n \geq 1} F(m - d, n)x^{m}y^{n}   \nonumber  \\
&\qquad  \qquad + \frac{1}{7}\sum_{m \geq 1} \sum_{n \geq 1} F(m, n - d_1)x^{m}y^{n}  +  \frac{1}{7} \sum_{m\geq 1}\sum_{n \geq 1} F(m, n - d_2)x^{m}y^{n}  \nonumber \\
&\qquad \qquad  + \frac{1}{7} \sum_{m \geq 1} \sum_{n \geq 1} F(m, n - d_1 - d_2)x^{m}y^{n} .
\end{align*}
By Lemma \ref{lemma: rewrite F(m -a, n - b)}, we can rewrite $A(x,y)$ as
\begin{align*}
A(x,y) \ = \  1 &+ \frac{1}{7}\Biggl[y^{d_1}\left( A(x,y) - 1 \right) + y^{d_2}\left(A(x,y) - 1\right) + y^{d_1 + d_2} \left(A(x,y) - 1\right)   \nonumber \\
&+  x^{d}\left(A(x,y) - 1\right) + \sum_{(D_1, D_2) \in \mathcal{D}} \left(\sum^{D_1}_{M = 1} \sum^{D_2}_{N = 1}x^M y^N + x^{D_1}y^{D_2} \left(A(x, y) - 1\right) \right)  \Biggr],
\end{align*}
where $\mathcal{D} = [d] \times [d_1, d_2, d_1 + d_2]$. We rearrange all terms containing $A(x,y)$ to the LHS and the rest of the terms to the RHS. We see the LHS equals
\[
A(x,y)\left(1 - \frac{1}{7} \left(x^{d}y^{d_1}+ x^{d}y^{d_2}+  x^{d}y^{d_1 + d_2} +  y^{d_1} +  y^{d_2} + y^{d_1 + d_2} +  x^{d} \right) \right),
\]
and RHS equals
\[
\left(1 - \frac{1}{7} \left(x^{d}y^{d_1}+ x^{d}y^{d_2}+  x^{d}y^{d_1 + d_2} +  y^{d_1} +  y^{d_2} + y^{d_1 + d_2} +  x^{d} \right) \right)  + \frac{1}{7} \sum_{(D_1, D_2) \in \mathcal{D}}\sum^{D_1}_{M = 1} \sum^{D_2}_{N = 1}x^M y^N.
\]

Let $\alpha = x^{d}$ and $\beta = y^{d_1} + y^{d_2} + y^{d_1 + d_2}$. We see that
\begin{align*}
A(x, y)\left(1- \frac{1}{7}\left(\alpha + \beta + \alpha \beta\right)\right) &\ = \  \left(1 - \frac{1}{7}\left(\alpha + \beta + \alpha \beta \right) \right) +  \frac{1}{7} \sum_{(D_1, D_2) \in \mathcal{D}}\sum^{D_1}_{M = 1} \sum^{D_2}_{N = 1}x^M y^N . 
\end{align*}
After solving  for $A(x, y)$, we  have 
\begin{equation*}
A(x,y) \ = \  1 + \frac{1}{7} B\left(\alpha, \beta, \frac{1}{7} \right)\sum_{(D_1, D_2) \in \mathcal{D}}\sum^{D_1}_{M = 1} \sum^{D_2}_{N = 1}x^M y^N,
\end{equation*}
where we recall that 
\[
B\left(\alpha, \beta, \frac{1}{7} \right) \ = \ \frac{1}{1 - \frac{1}{7}\left( \alpha + \beta + \alpha \beta\right)}.
\]

We first compute the coefficients of $B\left(\alpha, \beta, \frac{1}{7} \right)$. After a simple application of the geometric series formula and multinomial expansion, we see that
\begin{align*}
B\left(\alpha, \beta, \frac{1}{7} \right) &\ = \  \sum_{m, n \geq 0} \alpha^m \beta^n \sum^{\min(m, n)}_{\ell = 0} \binom{m + n - \ell}{m - \ell, n - \ell, \ell}\left(\frac{1}{7}\right)^{m + n - \ell} \nonumber \\
&\ = \ \sum_{m, n \geq 0}\alpha^m \beta^n D\left(m, n, \frac{1}{7} \right),
\end{align*}
where we recall that 
\begin{equation*}
D\left(m, n, \frac{1}{7} \right) = \sum_{\ell = 0}^{\min(m, n)}\binom{m + n - \ell}{m - \ell, n - \ell, \ell}\left(\frac{1}{7}\right)^{m + n - \ell}.
\end{equation*}

We now rewrite $B(\alpha,\beta,1/7)$ in terms of $x$ and $y$, where $\alpha = x^{d}$ and $\beta = y^{d_1} + y^{d_2} + y^{d_1 + d_2}$. After an application of multinomial expansion and a change of variables, we find
\begin{align}
B\left(\alpha, \beta, \frac{1}{7} \right) &\ = \ \sum_{m, n \geq 0}
x^{dm}\left(y^{d_1} + y^{d_2} + y^{d_1 + d_2} \right)^n D\left(m, n, \frac{1}{7} \right) \nonumber \\
&\ = \ \sum_{m, n \geq 0} x^{dm}  \sum_{\substack{i, j, k \geq 0 \\ i + j + k = n}} \binom{i + j + k}{i, j, k} y^{d_1(i +k) + d_2 (j + k)}   D\left(m, n, \frac{1}{7} \right).  
\label{eqs: method 2 denominator}
\end{align}
Now, we insert the expression for $B\left(\alpha,\beta, \frac{1}{7} \right)$ from Equation \eqref{eqs: method 2 denominator} into $A(x, y)$ and find that 
\begin{align*} 
A(x, y) &\ = \  1 +  \frac{1}{7} B\left(\alpha, \beta, \frac{1}{7} \right)\sum_{(D_1, D_2) \in \mathcal{D}}\sum^{D_1}_{M = 1} \sum^{D_2}_{N = 1}x^M y^N \\
&\ = \ 1 +  \frac{1}{7}\left(\sum_{m, n \geq 0} x^{dm}  \sum_{\substack{i, j, k \geq 0 \\ i + j + k = n}} \binom{i + j + k}{i, j, k} y^{d_1(i +k) + d_2 (j + k)}   D\left(m, n, \frac{1}{7} \right) \right) \\ 
&\qquad \qquad \cdot \left( \sum_{(D_1, D_2) \in \mathcal{D}} \sum_{M = 1}^{D_1} \sum^{D_2}_{N = 1} x^M y^N \right)\\
&\ = \ 
1 + \frac{1}{7}\Biggl( \sum_{m, n \geq 0} x^{dm + M} \sum_{(D_1, D_2) \in \mathcal{D}}\sum_{M = 1}^{D_1} \sum_{N = 1}^{D_2} \sum_{\substack{i, j, k \geq 0 \\ i + j + k = n}} \binom{i+ j + k}{i, j, k} \\
&\qquad \qquad \cdot y^{d_1(i + k) + d_2 (j + k) + N}  D\left(m, n, \frac{1}{7} \right) \Biggr).
\end{align*}
Let $I_1 = dm + M$ and $I_2 = d_1(i + k) + d_2 (j + k) + N$ for $1 \leq M \leq D_1$ and $1 \leq N \leq D_2$. We finally obtain the following generating function  
\begin{align*} 
A(x, y) &\ = \ 1 + \frac{1}{7}\Biggl(\sum_{(D_1, D_2) \in \mathcal{D}} \sum_{I_1, I_2 \geq 1}x^{I_1}y^{I_2} \\ 
&\qquad \qquad  \sum_{\substack{i, j, k \geq 0 \\ 1 \leq I_2 - d_1(i + k) - d_2(j + k) \leq D_2}}   \binom{i+ j + k}{i, j, k}D\left(\left\lfloor \frac{I_1 - 1}{d} \right\rfloor, i + j + k, \frac{1}{7}  \right) \Biggr). 
\end{align*}
Thus, we find the coefficients $[x^{I_1}y^{I_2}]A(x,y)$ as
\begin{equation*} 
\frac{1}{7} \sum_{(D_1, D_2) \in \mathcal{D}} \sum_{\substack{i, j, k \geq 0 \\ 1 \leq I_2 - d_1(i + k) - d_2(j + k) \leq D_2}}\binom{i+ j + k}{i, j, k}  
D\left(\left\lfloor \frac{I_1 - 1}{d} \right\rfloor, i + j + k, \frac{1}{7} \right).  
\end{equation*} 
\qed

\subsection{Proof of Theorem \ref{thm: finite sum for Game 3}}
We analyze Game 3.  Let $d_1, d_2 \geq 1$, so we have the following recurrence relation 
\begin{equation*}
F(m,n) =
\begin{cases}
\begin{aligned}
&\frac{1}{15}\biggl[ F(m -d_1-d_2, n -d_1 - d_2) \\
&\quad + F(m - d_1 - d_2, n - d_1) \\
&\quad + F(m - d_1 - d_2, n - d_2)  + F(m - d_1 - d_2, n) \\
&\quad + F(m - d_1, n - d_1 - d_2) + F(m - d_1, n - d_1) \\
&\quad + F(m - d_1, n - d_2) + F(m - d_1, n) \\ 
&\quad + F(m - d_2, n - d_1 - d_2) + F(m - d_2, n - d_1) \\
&\quad + F(m - d_2, n - d_2)   + F(m - d_2, n) \\ 
&\quad + F(m, n - d_1 - d_2) + F(m, n - d_1) \\
&\quad + F(m, n -d_2) \biggr]
\end{aligned}
& \text{ if } m,n \geq 1, \\
0 & \text{ if } m\leq0,n \geq 1  \\
& \text{ or } n\leq0,m \geq 1,  \\
1 & \text{ if } m, n \leq 0.
\end{cases}
\end{equation*}
Define the generating function $A(x, y) = \sum_{m, n \geq 0} F(m, n) x^m y^n$.  We have the following
\begin{align*}
A(x,y) &\ = \
1 + \frac{1}{15}
\Biggl[y^{d_1 + d_2}\left(A(x, y) - 1\right) + y^{d_1} \left(A(x, y) - 1\right) + y^{d_2} \left(A(x, y) - 1\right) \nonumber \\
&\qquad \qquad + x^{d_1 + d_2}\left( A(x,y) - 1\right) + x^{d_1} \left(A(x, y) - 1\right) + x^{d_2} \left( A(x,y) - 1\right) \nonumber \\
&\qquad \qquad +\sum_{(D_1, D_2) \in \mathcal{D}} \left(\sum^{D_1}_{M = 1} \sum^{D_2}_{N = 1}x^M y^N + x^{D_1}y^{D_2} \left(A(x, y) - 1\right) \right)  \Biggr],
\end{align*}
where $\mathcal{D} = [d_1, d_2, d_1 + d_2] \times [d_1, d_2, d_1 + d_2]$. We rearrange all terms containing $A(x,y)$ to the LHS and the rest of the terms to the RHS. We see the LHS equals
\begin{align*}
&A(x, y)\biggl(1 - \frac{1}{15}  \bigl( y^{d_1}  + y^{d_2}  +  y^{d_1 + d_2}  +  x^{d_1}  +  x^{d_2}  +  x^{d_1 + d_2} \nonumber \\
&\qquad\qquad +  x^{d_1}y^{d_1}  +  x^{d_1}y^{d_2}  + x^{d_1}y^{d_1 + d_2}  +  x^{d_2}y^{d_1}  + x^{d_2}y^{d_2}  +  x^{d_2}y^{d_1 + d_2}  \nonumber \\
&\qquad\qquad + x^{d_1 + d_2}y^{d_1}  +  x^{d_1 + d_2}y^{d_2} +  x^{d_1 + d_2}y^{d_1 + d_2} \bigr) \biggr),
\end{align*}
and RHS equals
\begin{align*}
& 1 - \frac{1}{15}  \bigl( y^{d_1}  + y^{d_2}  +  y^{d_1 + d_2}  +  x^{d_1}  +  x^{d_2}  +  x^{d_1 + d_2} \nonumber \\
&\qquad\qquad +  x^{d_1}y^{d_1}  +  x^{d_1}y^{d_2}  + x^{d_1}y^{d_1 + d_2}  +  x^{d_2}y^{d_1}  + x^{d_2}y^{d_2}  +  x^{d_2}y^{d_1 + d_2}  \nonumber \\
&\qquad\qquad + x^{d_1 + d_2}y^{d_1}  +  x^{d_1 + d_2}y^{d_2} +  x^{d_1 + d_2}y^{d_1 + d_2} \bigr) \nonumber \\
&\qquad \qquad + \frac{1}{15}\sum_{(D_1, D_2) \in \mathcal{D}}  \sum^{D_1}_{M = 1} \sum^{D_2}_{N = 1}x^M y^N.
\end{align*}

Let $\alpha = x^{d_1} + x^{d_2} + x^{d_1 + d_2}$ and $\beta = y^{d_1} + y^{d_2} + y^{d_1 + d_2}$. After solving for $A(x, y)$, we find
\begin{align}
A(x,y) &\ = \ 1 + \frac{\frac{1}{15} \sum_{(D_1, D_2) \in \mathcal{D}}  \sum^{D_1}_{M = 1} \sum^{D_2}_{N = 1}x^M y^N }{1 - \frac{1}{15}\left(\alpha + \beta + \alpha \beta \right)} \nonumber \\
&\ = \ 1 + \frac{1}{15}B\left(\alpha, \beta, \frac{1}{15}\right)\sum_{(D_1, D_2) \in \mathcal{D}}  \sum^{D_1}_{M = 1} \sum^{D_2}_{N = 1}x^M y^N,
\label{eqs: method 3 A(x,y) in terms of B}
\end{align}
where
\begin{equation*}
B\left(\alpha, \beta, \frac{1}{15}\right) \ = \  \frac{1}{1 - \frac{1}{15}\left(\alpha + \beta + \alpha \beta \right)}.
\end{equation*}

We now expand  $B(\alpha, \beta, \frac{1}{15})$ and write its coefficients in terms of $x$ and $y$. After a simple application of the geometric series formula and multinomial expansion, we first see that
\begin{equation}
B\left(\alpha, \beta, \frac{1}{15}\right) \ = \ \sum_{m, n \geq 0} \alpha^m \beta^n D\left(m, n, \frac{1}{15} \right),
\label{eqs: method 3 geometric series}
\end{equation}
where
\begin{equation*}
D\left(m, n, \frac{1}{15} \right) \ = \ \sum^{\min(m, n)}_{\ell = 0} \binom{m + n - \ell}{m - \ell, n - \ell, \ell}\left(\frac{1}{15}\right)^{m + n - \ell}.
\end{equation*}
Substituting $\alpha = x^{d_1} + x^{d_2} + x^{d_1 + d_2}$ and $\beta = y^{d_1} + y^{d_2} + y^{d_1 + d_2}$ into Equation \eqref{eqs: method 3 geometric series}, we arrive at
\begin{align}
B\left(\alpha, \beta, \frac{1}{15} \right)
&= \sum_{m, n \geq 0}
   \left(x^{d_1} + x^{d_2} + x^{d_1 + d_2}\right)^m
   \left(y^{d_1} + y^{d_2} + y^{d_1 + d_2}\right)^n \nonumber \\ 
   &\cdot 
   D\left(m, n, \frac{1}{15} \right)
   \nonumber \\
&= \sum_{m, n \geq 0} \Biggl[
   \sum_{\substack{r,s,t \geq 0 \\ r + s + t = m}}
   \binom{r + s + t}{r,s,t}
   x^{d_1(r + t) + d_2(s + t)}
   \nonumber \\
&\qquad
   \sum_{\substack{i, j, k \geq 0 \\ i + j + k = n}}
   \binom{i + j + k}{i, j, k}
   y^{d_1(i + k) + d_2(j + k)}
   D\left(r + s + t, i + j + k, \frac{1}{15} \right)\Biggr] . 
\label{eqs: method 3 denominator} 
\end{align}
We have now re-parametrized $B(\alpha, \beta, \frac{1}{15})$ in terms of $\alpha = x^{d_1} + x^{d_2} + x^{d_1 + d_2}$ and $\beta = y^{d_1} + d^{d_2} + y^{d_1 + d_2}$. We insert $B(\alpha, \beta, \frac{1}{15})$ from Equation \eqref{eqs: method 3 denominator} into Equation \eqref{eqs: method 3 A(x,y) in terms of B} and find that  
\begin{align*} 
A(x, y) &\ = \ 1 +  \frac{1}{15} B\left(\alpha, \beta, \frac{1}{15}\right) \cdot \left( \sum_{(D_1, D_2) \in \mathcal{D}} \sum_{M = 1}^{D_1} \sum_{N = 1}^{D_2}x^M y^N \right) \nonumber \\ 
&\ = \ 1 + \frac{1}{15}\Biggl(   \sum_{m, n \geq 0} \Biggl[
   \sum_{\substack{r,s,t \geq 0 \\ r + s + t = m}}
   \binom{r + s + t}{r,s,t}
   x^{d_1(r + t) + d_2(s + t)}
   \nonumber \\
&\qquad
   \sum_{\substack{i, j, k \geq 0 \\ i + j + k = n}}
   \binom{i + j + k}{i, j, k}
   y^{d_1(i + k) + d_2(j + k)}
   D\left(r + s + t, i + j + k, \frac{1}{15} \right)\Biggr]\Biggr) \\ 
   &\qquad \left( \sum_{(D_1, D_2) \in \mathcal{D}} \sum_{M = 1}^{D_1} \sum_{N = 1}^{D_2}x^M y^N \right) \nonumber \\ 
&\ = \ 1 + \frac{1}{15}\Biggl(\sum_{(D_1, D_2) \in \mathcal{D}}\sum_{M = 1}^{D_1}\sum_{N = 1}^{D_2} \sum_{m, n \geq 0} \sum_{\substack{r, s, t \geq 0 \\ r + s + t = m}} \binom{r + s + t}{r,s,t} x^{d_1(r + t) + d_2 (s + t) + M} \nonumber \\ 
&\qquad \qquad \cdot \sum_{\substack{i, j, k \geq 0 \\ i + j + k = n}}\binom{i + j + k}{i, j, k} y^{d_1(i + k) + d_2(j + k) + N} D\left(r + s + t, i + j + k, \frac{1}{15} \right)\Biggl) \nonumber \\
&\ = \ 1 + \frac{1}{15}\Biggl(\sum_{(D_1, D_2) \in \mathcal{D}} \sum_{I_1, I_2 \geq 1}  x^{I_1}y^{I_2} \sum_{\substack{r, s, t \geq 0 \\ 1 \leq I_1 - d_1(r + t) - d_2(s + t) \leq D_1}} \binom{r + s + t}{r,s,t}  \nonumber \\
&\qquad \qquad \qquad \sum_{\substack{i, j, k \geq 0 \\ 1 \leq I_2 - d_1(i + k) - d_2(j + k) \leq D_2}}\binom{i + j + k}{i, j, k}   D\left(r + s + t, i + j + k, \frac{1}{15} \right)\Biggl) . 
\end{align*}
Let $I_1 = d_1(r + t) + d_2(s + t) + M$ and $I_2 = d_1(i + k) + d_2(j + k) + N$ for $1\leq M \leq D_1$ and $1 \leq N \leq D_2$. Thus, the coefficients $[x^{I_1}y^{I_2}]A(x, y)$ are 
\begin{align*} 
&\frac{1}{15}\Biggl(\sum_{(D_1, D_2) \in \mathcal{D}}  \sum_{\substack{r, s, t \geq 0 \\ 1 \leq I_1 - d_1(r + t) - d_2(s + t) \leq D_1}} \sum_{\substack{i, j, k \geq 0 \\ 1 \leq I_2 - d_1(i + k) - d_2(j + k) \leq D_2}} \binom{r + s + t}{r,s,t}  \nonumber \\
&\qquad \qquad \qquad \cdot  \binom{i + j + k}{i, j, k}   D\left(r + s + t, i + j + k, \frac{1}{15} \right)\Biggl) . 
\end{align*}
\qed

\section*{Acknowledgments} The generalized M\&M Game was created by the second and sixth named authors at the Zassenhaus Groups and Friends Conference at Texas State University in June 2024; we thank the organizers for fostering an environment conducive to exploration. This work began in the 2024 Polymath Junior Program, supported in part by NSF Grant DMS2341670. We thank the students and colleagues at Polymath Junior, the participants in the Research Special Session at the 2025 Joint Mathematics Meetings, and the anonymous referee for helpful comments. We also thank Evan Li and Andrew Mou for their help with earlier versions of this paper. We used Anthropic's freely available Claude Sonnet 4.6 to check for grammatical errors and unclear phrasing, and we revised the text in several places based on its suggestions.

{\footnotesize  
\medskip
\medskip
\vspace*{1mm} 
 
\noindent {\it Snehesh Das}\\  
University of California, Los Angeles\\
405 Hilgard Avenue \\
Los Angeles, CA 90024 \\
E-mail: {\tt snehesh2016@g.ucla.edu}\\ \\  

\noindent {\it Steven J. Miller}\\  
Williams College \\
880 Main Street\\
Williamstown, MA 01267\\
E-mail: {\tt sjm1@williams.edu}\\ \\

\noindent {\it Geremias Polanco}\\  
Smith College \\
10 Elm Street \\
Northampton, MA 01063 \\
E-mail: {\tt gpolanco@smith.edu}\\ \\

\noindent {\it Yilong Wu}\\  
Williams College \\
880 Main Street\\
Williamstown, MA 01267 \\
E-mail: {\tt mw17@williams.edu}\\ \\

\noindent {\it Xiaochen Wang}\\  
Hong Kong University of Science and Technology \\
Hong Kong \\
E-mail: {\tt xwanggk@connect.ust.hk} \\ \\

\noindent {\it  April Yang
}\\  
Doris Miller Middle School  \\
301 Foxtail Run \\
San Marcos, TX 78666 \\ 
E-mail: {\tt apr.yan077@gae.smcisd.net} \\ \\   
\noindent {\it  Chris Yao
}\\  
University of California, Berkeley  \\
110 Sproul Hall \\
Berkeley, CA 94720 \\ 
E-mail: {\tt chris.yao@berkeley.edu} \\ \\


}

\vspace*{1mm}\noindent\footnotesize{\date{ {\bf Received}: January 31, 2025\;\;\;{\bf Accepted}: July 2, 2026}}\\
\vspace*{1mm}\noindent\footnotesize{\date{  {\bf Communicated by Serban Raianu}}}


{\footnotesize
\begin{thebibliography}{00}

\bibitem{badinski2017m}
I. Badinski, C. Huffaker, N. McCue, C.N. Miller, K.S. Miller, S.J. Miller, M. Stone,
The M\&M Game: From Morsels to Modern Mathematics,
{\it Math. Mag.}, {\bf 90} (2017), 197--207.

\bibitem{Banderier_2005}
C. Banderier, S. Schwer,
Why Delannoy numbers?,
{\it J. Statist. Plann. Inference}, {\bf 135} (2005), 40--54.

\bibitem{dhar2018comparison}
M. Dhar, P. Bhattacharya,
Comparison of the logistic and the Gompertz curve under different constraints,
{\it JSMS}, {\bf 21} (2018), 1189--1210.

\bibitem{marquardt1963algorithm}
D.W. Marquardt,
An algorithm for least-squares estimation of nonlinear parameters,
{\it SIAP}, {\bf 11} (1963), 431--441.

\bibitem{mcbook}
A.B. Owen,
{\it Monte Carlo Theory, Methods and Examples},
available online at the URL:
\url{https://artowen.su.domains/mc/}. 

\bibitem{rice2006mathematical}
J.A. Rice,
{\it Mathematical Statistics and Data Analysis},
3rd ed., Cengage Learning, 2006.

\bibitem{wilf_generatingfunctionology}
H.S. Wilf,
{\it Generatingfunctionology},
3rd ed., A K Peters, 2005.
\end{thebibliography}
}
\end{document}